\documentclass[smallextended]{svjour3}       % onecolumn (second format)
\smartqed  % flush right qed marks, e.g. at end of proof
\usepackage{graphicx}
\usepackage{amssymb,amsfonts}
\usepackage[all,knot]{xy}
\usepackage{enumerate}
\usepackage{mathrsfs}
\usepackage{braket}

\journalname{Quantum Information Processing}
\begin{document}

\title{Topological Aspects of Quantum Entanglement}


\author{Louis H. Kauffman         \and
        Eshan Mehrotra %etc.
}


\institute{L. Kauffman \at
              Mathematics Department, UIC, 851 South Morgan Street, Chicago IL, 60607-7045 \\
              Tel.: +312-996-3066\\
              Fax: +3129961491\\
              \email{kauffman@uic.edu}\\
              and\\
              Department of Mechanics and Mathematics\\ Novosibirsk State University\\
              Novosibirsk, Russia
                         %  \\
           \and
           E. Mehrotra \at
              Illinois Mathematics and Science Academy, 1500 Sullivan Road, Aurora IL, 60506 \\
							\email{emehrotra@imsa.edu}
}

\date{Received: date / Accepted: date}

\maketitle

\begin{abstract}
Kauffman and Lomonaco in \cite{notsoshort} and \cite{May} explored the idea of understanding quantum entanglement (the non-local correlation of certain properties of particles) topologically by viewing unitary entangling operators as braiding operators. In \cite{amsshort}, it is shown that entanglement is a necessary condition for forming non-trivial invariants of knots from braid closures via solutions to the Yang-Baxter Equation. We show that the arguments used by \cite{amsshort} generalize to essentially the same results for quantum invariant state summation models of knots. In one case (the unoriented swap case) we give an example of a Yang-Baxter operator, and associated quantum invariant, that can detect the Hopf link. Again this is analogous to the results of \cite{amsshort}. We also give a class of $R$ matrices that are entangling and are weak invariants of classical knots and links yet strong invariants of virtual knots and links. We also give an example of an $SU(2)$ representation of the three-strand braid group that models the Jones polynomial for closures of three-strand braids. This invariant is a quantum model for the Jones polynomial restricted to three strand braids, and it does not involve quantum entanglement. These relationships between topological braiding and quantum entanglement can be used as a framework for future work in understanding the properties of entangling gates in topological quantum computing.  The paper ends with a discussion of the Aravind hypothesis about the direct relationship of knots and quantum entanglement, and the $ER = EPR$ hypothesis about the 
relationship of quantum entanglement with the connectivity of space. We describe how, given a background space and a quantum tensor network, to construct a new topological space, that welds the network and the background space together. This construction embodies the principle that quantum entanglement and topological connectivity are 
intimately related.

\keywords{topological entanglement \and quantum entanglement \and Yang-Baxter operator \and state summation \and quantum link invariant}
\end{abstract}

\section{Introduction}
\label{intro}
The purpose of this paper is to explore several phenomena that relate topology and quantum entanglement. Braiding operators are topological objects, while unitary operators are primarily used in the realm of quantum mechanics. This paper establishes a relationship between the two. We first examinine a quantum gate $R$ which is both entangling and unitary. Such gates are useful for quantum computation. Second, we choose an $R$ that satisfies the Yang-Baxter equation and determine the relation between entangling $R$'s and detecting knotting and linking. We show in this paper that non-entangling Yang-Baxter operators cannot  form non-trivial invariants of knots in the oriented and unoriented cases of quantum state summations. There do exist cases where we can construct non-trivial invariants of knots and links from unitary transformations where the operators are not entangling. For example, the Jones polynomial {\cite{Jones,Jones1,Jones2,Jones3,Jones4} for three strand braids can be extracted from computations that involve only a single qubit \cite{quantumcompute}. See Section 5 of the present paper.\\

Section 2 of this paper explicates the relationship between unitary operators and braiding operators, while also providing a brief introduction to the theory of quantum link invariants. Section 3 shows that the results of \cite{amsshort} generalize to unoriented quantum invariant state summations in the so-called product case. In the swap case, considered in \cite{amsshort}, the Markov trace method for constructing the proposed link invariant does not generalize to a quantum summation of the kind we consider but we nevertheless give an example of a Yang-Baxter operator in this case that can detect the Hopf Link. This lack of correspondence is interesting in its own right, and is discussed in this section.  Section 4 shows that non-trivial invariants with non-entangling Yang-Baxter operators cannot be constructed in the oriented case. Section 5  describes how the Jones polynomial can still arise in systems that lack quantum entanglement. Section 6 describes how unitary $R$ matrix solutions to the bracket state summation are unentangling. Section 7 establishes a potential relationship between quantum entanglement and virtual knots and links. Section 8 is wider discussion of the relationship between topology and entanglment. We discuss the Aravind hypothesis that suggests that knots and links themselves may be connected more directly with quantum entanglement, and we discuss the $ER = EPR$ hypothesis of Leonard Susskind and his collaborators that suggests that the connectivity of space itself is directly related to quantum entanglement. We illustrate this ideas of connectivity by showing how the tensor networks for entangled states (in the sense of the networks used in the present paper) can be used to both indicate this new connectivity and can be welded to the given space by adding points for the entangled states and new neighborhoods to extend the topology. We describe how, given a background space and a quantum tensor network, to construct a new topological space, that welds the network and the background space together. This part of the paper is intended to be brief and will be expanded further in subsequent work. Finally,  Section 9 concludes the paper with a discussion of the ideas and concepts that have arisen during the course of this research.  The Appendix proves an 
important Lemma for our analysis of link invariants in the earlier parts of the paper.\\

\section{Characteristics of Unitary Operators and the Artin Braid Group}

We begin by describing the Artin braid group \cite{Birman}. Figure~\ref{fig:operator} shows the elements of this group. An $n$-stranded braid is a collection of $n$ strings extending from one row of $n$ points to another row of $n$ points,  with each cross section of the braid consisting of $n$ points. The $n$-strand braid group $B_n$ is generated by $\sigma_1,... , \sigma_{n-1}$ where $\sigma_i$ is a twist of the $i$ and $i+1$ strands as shown in Figure~\ref{fig:operator}. The relations on these generators are given by $\sigma_i \sigma_j = \sigma_j \sigma_i$ for $|i-j|>1$ and $\sigma_i \sigma_{i+1} \sigma_i = \sigma_{i+1} \sigma_i \sigma_{i+1}$ for $i=1,..., n-2$. Braid multiplication is defined by attaching the initial points of one braid to the end points of the other. Under topological equivalence, this multiplication operation gives the Artin braid group $B_{n}$ for $n$-stranded braids. Figure~\ref{fig:inverse} shows two 2-strand braids and a respective braid multiplication between them that demonstrates multiplicative inverse. \\

\begin{figure}
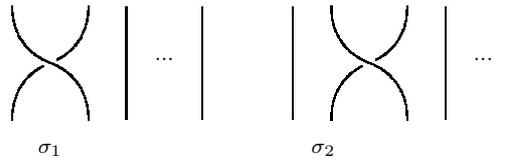

	\[
	\xy
	\vcross~{(-5,5)}{(5,5)}{(-5,-10)}{(5,-10)}; 
	(10,5)*{}; (10,-10)*{} **\dir{-};
	(20,5)*{}; (20,-10)*{} **\dir{-};
	(15, -2)*{. . .};
	(0,-14)*{\sigma_{1}};
	\endxy
	\qquad \qquad
	\xy
	(-10,5)*{}; (-10,-10)*{} **\dir{-};
	\vcross~{(-5,5)}{(5,5)}{(-5,-10)}{(5,-10)}; 
	(10,5)*{}; (10,-10)*{} **\dir{-};
	(20,5)*{}; (20,-10)*{} **\dir{-};
	(15, -2)*{. . .};
	(-6,-14)*{\sigma_{2}};
	\endxy
	\]
	\caption{The n-stranded braiding operators.}
	\label{fig:operator}
\end{figure}

\begin{figure}
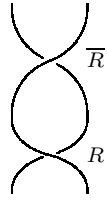


	\[ \xy 
	\vtwist~{(-5,5)}{(5,5)}{(-5,-10)}{(5,-10)}; 
	\vcross~{(-5,-10)}{(5,-10)}{(-5,-20)}{(5,-20)};
	(6,-2)*{\overline{R}};
	(6,-15)*{R};
	\endxy \]
	\caption{Two-strand braid inverses.}
	\label{fig:inverse}
\end{figure}

We can study quantum entanglement and topological quantum information by examining unitary representations of the Artin braid group. In such a representation each braid is mapped to a unitary operator. 
Given such a representation, we can examine the entangling capacity of the braiding operators.  That is, we can calculate whether they can take unentangled states to entangled states. It is also possible to use such a braiding representation to create topological invariants of knots, links and braids. Thus one can, in principle, compare the power of such a representation to detect knots and links with the quantum entangling capacity of the operators in the representation.\\

Consider representations of the braid group such that for a single twist, as in the lower half of Figure~\ref{fig:inverse}, there is an associated operator

$$
	R: V \otimes V \rightarrow V \otimes V.
$$

\noindent In the above operator, V is a complex vector space (In this case we take $V$ to be two dimensional so that it can hold a single qubit of information. In general the restriction is not necessary.). The two input and two output lines in the braid (see $R$ in Figure~\ref{fig:cupscaps}) are representative of the fact that the operator $R$ is defined on the tensor product of complex vector spaces. Thus, the top endpoints of $R$ as shown in Figure~\ref{fig:cupscaps} represent $V \otimes V$ as the domain of $R$, and the bottom endpoints of $R$ represent $V \otimes V$ as the range of $R$. The diagram in Figure~\ref{fig:yangbaxt} shows mappings of $V \otimes V \otimes V$ to itself. This relation is the Yang-Baxter equation \cite{Baxter}. Algebraically with $I$ representing the identity on $V$, the equation reads as follows:

\begin{figure}
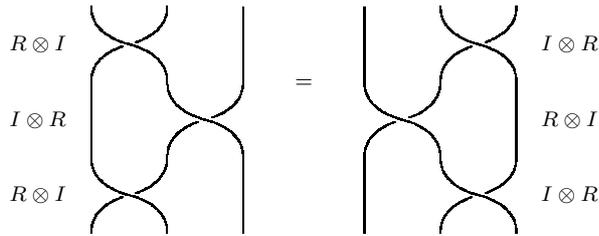

		
	\[ 
	\xy 
	\vcross~{(-5,10)}{(5,10)}{(-5,0)}{(5,0)}; 
	(15,10)*{}; (15,0)*{} **\dir{-};
	\vcross~{(5,0)}{(15,0)}{(5,-10)}{(15,-10)}; 
	(-5,0)*{}; (-5,-10)*{} **\dir{-}; 
	\vcross~{(-5,-10)}{(5,-10)}{(-5,-20)}{(5,-20)}; 
	(15,-10)*{}; (15,-20)*{} **\dir{-}; 
	(-12, 5)*{R \otimes I};
	(-12, -5)*{I \otimes R};	
	(-12,-15)*{R \otimes I};
	\endxy
	\qquad = \qquad
	\xy 
	\vcross~{(5,10)}{(15,10)}{(5,0)}{(15,0)}; 
	(-5,10)*{}; (-5,0)*{} **\dir{-};
	\vcross~{(-5,0)}{(5,0)}{(-5,-10)}{(5,-10)}; 
	(15,0)*{}; (15,-10)*{} **\dir{-}; 
	\vcross~{(5,-10)}{(15,-10)}{(5,-20)}{(15,-20)}; 
	(-5,-10)*{}; (-5,-20)*{} **\dir{-}; 
	(22, 5)*{I \otimes R};
	(22, -5)*{R \otimes I};	
	(22,-15)*{I \otimes R};
	\endxy
	\]
	\caption{The Yang-Baxter equation.}
	\label{fig:yangbaxt}
\end{figure}

$$
	(R \otimes I)(I \otimes R)(R \otimes I) = (I \otimes R)(R \otimes I)(I \otimes R).
$$

\noindent This equation represents the fundamental topological relation in the Artin braid group. If $R$ satisfies the Yang-Baxter equation and is invertible, then we can define a representation $\tau$ of the braid group by 

$$
	\tau(\sigma_{k}) = I \otimes ... \otimes I \otimes R \otimes I ... \otimes I,
$$

\noindent where $R$ occupies the $k$ and $k+1$ places in the above tensor product. If $R$ is unitary, then this is a unitary representation of the braid group. 
Since the basic operator $R$ operates on $V \otimes V$, a tensor product of qubit spaces, it is possible to measure whether it is an entangling operator. In previous work \cite{notsoshort} we found that there appears to be a relationship between such entangling capacity and the ability to use $R$ to produce a non-trivial invariant of knots and links. Alagic, Jarret and Jordan \cite{amsshort} proved, using Markov trace models \cite{Birman} for link invariants associated with braids, that if the operator $R$ is not an entangling operator, then the corresponding knot invariants are trivial. In this paper, we corroborate their results for state sum models (defined on general link diagrams).\\

It should be remarked that what we have above called Markov trace models for link invariants are based on a fundamental theorem of J. W. Alexander \cite{Alex} that states that any knot or link has a representation as the closure of a braid.  A braid, as depicted above, can be {\it closed} by attaching the upper strands to the lower strands by a parallel bundle of non-crossing strands that is positioned next to the given braid. The result of the closure is that the diagram of the closed braid has the appearance of a bundle of strands that proceeds circularly around an axis perpendicular to the plane. Alexander shows how to isotope any knot of link into such a form.
It is then the case that a given link can be obtained as the closure of different braids. The Markov Theorem \cite{Birman} gives an equivalence relation on braids so that two braids close to the same knot or link if and only if they are 
Markov equivalent. By constructing functions on braids that are invariant under the generating moves for Markov equivalence, one produces Markov trace invariants of knots and links. Such invariants can be constructed from solutions $R$ to the Yang-Baxter equation and some extra information. This approach is used by Alagic, Jarret and Jordan \cite{amsshort}. \\

In the next section, we describe quantum link invariants and prove theorems showing their limitations when built with non-entangling solutions to the Yang-Baxter equation. The class of quantum link invariant state sum models is very closely related to Markov trace models, but one does not need to transform the knot or link to a closed braid form.\\ 

\subsection{\bf Quantum Link Invariants}

We now describe how invariants of knots and links can be constructed by arranging knots and links with respect to a given direction in the plane denoted as \textit{time}. Consider the circle in a spacetime plane with time on the vertical axis and space on the horizontal axis. This is shown in Figure~\ref{fig:spacetime}. The circle, under this paradigm, represents a vacuum to vacuum process that depicts the creation of two \textit{particles}  and their subsequent annihilation. The two parts of this process are represented by a creation \textit{cup} (the bottom half of the circle) and an annihilation \textit{cap} (the top half of the circle). We can then consider the amplitude of this process given by $\braket{cap|cup}$. Since the diagram for the creation of the two particles ends in two separate points, it is natural to take a vector space of the form $V \otimes V$ as the target for the bra and as the domain of the ket. We imagine at least one particle property being catalogued by each factor of the tensor. We use this physical metaphor to describe the model. It is understood that the model applies to mathematical or topological situations where time is just a convenient parameter and particles are just matrix indices.  Knot and link invariants built in this framework are called {\it quantum link invariants}
because the numerical value of the invariant can be interpreted as a (generalized) amplitude for the vacuum to vacuum process represented by the link diagram. We give the details of this formulation below.\\

\begin{figure}
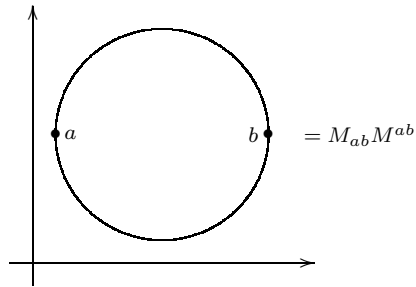


	\[
	\xy
	(-14,0)*+{\bullet}="4"; 
	(14,0)*+{\bullet}="6";
	(0,0)*\xycircle(14,14){-};
	{\ar (-20,-17)*{}; (20,-17)*{}};
	{\ar (-17,-20)*{}; (-17,17)*{}};
	(-12,0)*{a};
	(12,0)*{b};
	(26,0)*{= M_{ab}M^{ab}};
	\endxy 
	\]
	\caption{The quantum link invariant based evaluation of a circle in spacetime.}
	\label{fig:spacetime}
	
\end{figure}

We shall call a link diagram arranged with respect to a direction in time a {\it Morse diagram}. Note that, generically, in a Morse diagram, a horizontal line in the plane intersects the diagram transversely in a finite collection of 
points. Special points or {\it critical points} consist in maxima and minima in the diagram, and the places where a crossing appears in the diagram. We can transform any link diagram into a Morse diagram by an isotopy of the plane and so all knots and links are represented by Morse diagrams. Before going further with Morse diagrams, we first recall that two diagrams, regarded as projections of knots or links in three-space, are equivalent by Reidemeister moves as shown in Figure~\ref{reid}. This result, due to Reidemeister, Alexander and Briggs \cite{Reid}, implies that the equivalence classes of diagrams generated by the Reidemeister moves classify the topological types of knots and links in three-dimensional space. In order to work with Morse diagrams, we use a refomulation of the Reidemeister Theorem that utilizes the move types shown in 
Figure~\ref{morsemoves}. The reformulation of the Reidemeister theorem \cite{RT1,RT2,TM,Yetter} states that two Morse link diagrams are equivalent via the {\it Morse moves} of Figure~\ref{morsemoves} if and only if they are regularly isotopic. A good reference for the details of this theorem based on Reidemeister's orginal approach can be found in the paper by David Yetter \cite{Yetter}. {\it Regular isotopy} is the equivalence relation on diagrams generated by the second and third Reidemeister moves. Thus Morse diagrams and their moves give a complete formalism for the regular isotopy classification of standard knot and link diagrams. Regular isotopy invariance is often the most convenient method for studying knots and links. Invariants of regular isotopy can often be normalized to produce invariants of {\it ambient isotopy} (the equivalence relation generated by all three Reidemeister moves). In the following we shall detail how to use solutions of the Yang-Baxter equation to produce invariants of regular isotopy for Morse diagrams.\\

The strategy for this method to produce invariants is illustrated in Figure~\ref{jordanamp} and Figure~\ref{morseamp}. In the following we explain the use of Morse diagrams for producing link invariants. The original approach, due to Reshetikhin and Turaev \cite{RT1,RT2}, is formulated using the oriented tangle category. Our approach describes the analogous structure for unoriented diagrams and can be used as well for oriented diagrams. We divide the Morse diagram into parts that are the shape of a maxima, a minima or a crossing.
We associate matrices $M^{ab}$ to minima, $M_{ab}$ to maxima and $R^{ab}_{cd}$ to crossings. Each choice of indices for any matrix gives a scalar quantity for the corresponding matrix entry. The diagram yields,
as in  Figure~\ref{morseamp}, a product of these scalars with every index repeated twice. One then takes the summation of these products over all choices of indices. The resulting state summation $Z_{K}$ is the quantum link amplitude. In our physical metaphor, this is the quantum amplitude for the vacuum to vacuum process the involves the creation of particles via minima, the interaction of particles at the crossings and annihilations of particles at the maxima. The matrices must satisfy a collection of equations that correspond to the moves on Morse diagrams. We detail these equations and the correspondences below.\\

\begin{figure}
     \begin{center}
     \begin{tabular}{c}
     \includegraphics[width=6cm]{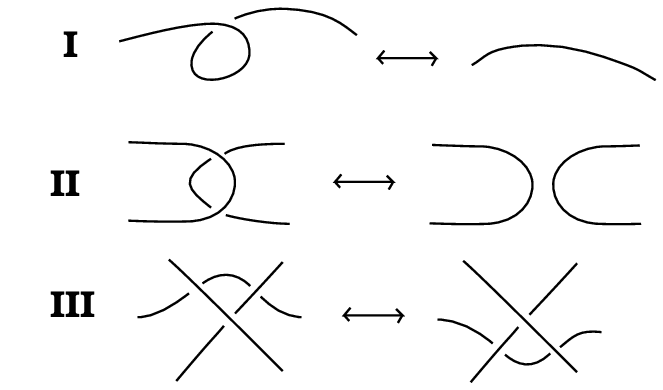}
     \end{tabular}
     \caption{\bf Classical Reidemeister Moves}
     \label{reid}
\end{center}
\end{figure}

\begin{figure}
     \begin{center}
     \begin{tabular}{c}
     \includegraphics[width=6cm]{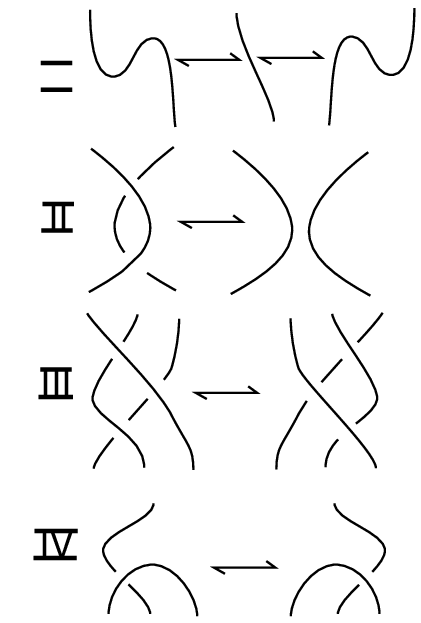}
     \end{tabular}
     \caption{\bf Regular Isotopy With Respect to a Vertical Direction}
     \label{morsemoves}
\end{center}
\end{figure}

\begin{figure}
     \begin{center}
     \begin{tabular}{c}
     \includegraphics[width=6cm]{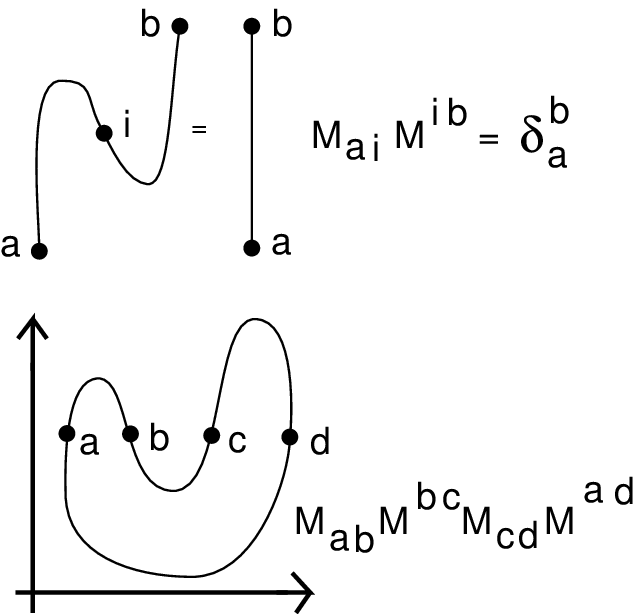}
     \end{tabular}
     \caption{\bf Jordan Curve Amplitude}
     \label{jordanamp}
\end{center}
\end{figure}

\begin{figure}
     \begin{center}
     \begin{tabular}{c}
     \includegraphics[width=6cm]{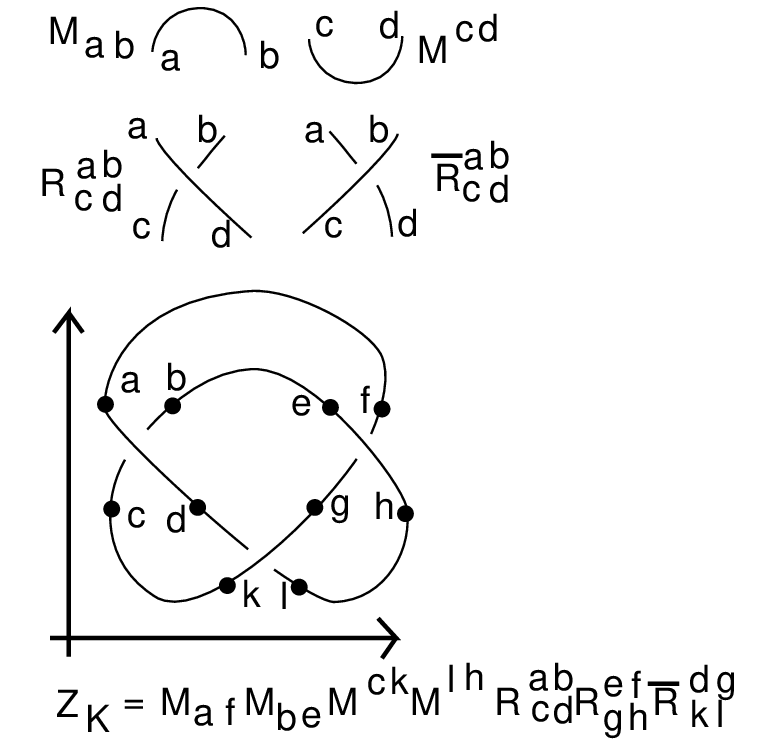}
     \end{tabular}
     \caption{\bf Amplitude for a Morse Diagram}
     \label{morseamp}
\end{center}
\end{figure}

All crossings in a link diagram are represented by transversal intersections. Any non-self-intersecting differentiable curve (for embedded curves and for transversely intersecting immersed curves) can be rigidly rotated until it is in general position with respect to the vertical. A curve without intersections  is then seen to decompose into an interconnection of minima and maxima. We can evaluate an amplitude for any curve in this general position with respect to a vertical direction. Any simple closed curve in the plane is isotopic to a circle, by the Jordan Curve Theorem. If these are topological amplitudes, then the value for any simple closed curve should be equal to the amplitude of the circle. In order to find conditions for the creation and annihilation operators that ensure amplitudes that respect topological equivalence, isotopies of simple closed curves are generated by the cancellation of adjacent maxima and minima. Specifically, let ${e_1, e_2, ..., e_n}$ be a basis for $V$. Let $e_{ab} = e_a \otimes e_b$ denote the elements of the tensor basis for $V \otimes V$. Then, there are matrices $M_{ab}$ and $M^{ab}$ such that 

\begin{figure}
	\[ 
	\xy 
	(0,0)*{}="A"; 
	(10,0)*{}="B"; 
	"A"; "B" **\crv{(5,5)}; 
	(18,0)*{}="C"; 
	(28,0)*{}="D"; 
	"C"; "D" **\crv{(23,-5)}; 
	\vcross~{(36,3)}{(46,3)}{(36,-3)}{(46,-3)};
	\vtwist~{(54,3)}{(64,3)}{(54,-3)}{(64,-3)};
	(5,-12)*{Cap};
	(23,-12)*{Cup};
	(41,-12)*{R};
	(59,-12)*{\overline{R}};
	\endxy
	\]
	\caption{}
	\label{fig:cupscaps}
\end{figure}

$$
\ket{cup}(1) = \sum{M^{ab}e_{ab}},
$$

\noindent with the summation taken over all values of $a$ and $b$ from $1$ to $n$. Similarly, $\bra{cap}$ is described by 

$$
\bra{cap}(e_{ab}) = M_{ab}.
$$

\noindent Thus the amplitude for the circle is 

$$
\braket{cap|cup}(1) = \bra{cap} \sum{M^{ab}e_{ab}} = \sum{M^{ab} \bra{cap}(e_{ab})} = \sum{M^{ab} M_{ab}}.
$$

In general, the value of the amplitude on a simple closed curve is obtained by translating it into an ``abstract tensor expression" using $M^{ab}$ and $M_{ab}$, and then summing over the products for all cases of repeated indices. Note that here the value ``1" corresponds to the vacuum. For example in Figure~\ref{jordanamp} we write down a more complex amplitude for a Jordan curve in the lower part of the figure. We also illustrate
a topological relation on the matrices that will ensure that this evaluation is the same as the circle evaluation above. This topological relation is just that the matrices $M^{ab}$ and $M_{cd}$ are inverses in the sense that 

$$
\sum_{i} M_{ai}M^{ib} = \delta^{b}_{a},
$$

\noindent where $\delta^{b}_{a}$ denotes the identity matrix. This equation is illustrated diagrammatically in Figure~\ref{jordanamp}. \\

One of our simplest choices is to take a $2 \times 2$ matrix $M$ such that $M^{2} = I$, where $I$ is the identity matrix. Then the entries of $M$ can be used for both the cup and the cap. The value for a loop is then equal to the sum of the squares of the entries of $M$:

$$
\braket{cap|cup} = \sum M^{ab} M_{ab} = \sum M_{ab} M_{ab} = \sum M^{2}_{ab}.
$$

Any knot or link can be represented by a picture that is configured with respect to a vertical direction in the plane. The picture decomposes into minima (creations), maxima (annihilations), and crossings of the two types shown in Figure~\ref{morseamp} and Figure~\ref{fig:cupscaps}. Here the knots and links are unoriented. Any knot or link can be written as a composition of these fragments, and consequently a choice of such mappings determines an amplitude for knots and links. In order for such an amplitude to be topological (i.e. an invariant of regular isotopy of the equivalence relation generated by the second and third classical Reidemeister moves) we want it to be invariant under a list of local moves as shown in  Figure~\ref{fig:six}, Figure~\ref{fig:seven}, Figure~\ref{fig:eight}, Figure~\ref{fig:nine}. \\

We now give an explanation of the algebraic and topological equations shown in these figures. Figure~\ref{fig:six} is the cancellation of maxima and minima. Figure~\ref{fig:seven} corresponds to the second Reidemeister move. Figure~\ref{fig:eight} is the Yang-Baxter equation. Figure~\ref{fig:nine} demonstrates that a line can move across a minimum (similar equations can be formulated for a line moving across a maximum). In each figure we have given the corresponding equation for the cup, cap and crossing matrix elements. If these equations are taken purely abstractly then they indicate a necessary and sufficient condition for a state sum of this type to be  
an invariant of regular isotopy. In order to produce an invariant, it is sufficient that the matrices satisfy these conditions. Such an invariant is not necessarily a complete invariant of regular isotopy, and to this date no one has 
produced such a complete invariant other than the formalism itself.\\

\begin{figure}
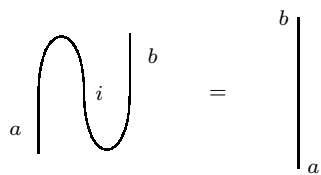

	\[
	\xy (-6,-8)*{}="1E"; (-6,0)*{}="1"; (0,0)*{}="2"; (6,0)*{}="3"; 
	(6,8)*{}="3B"; 
	"2";"1" **\crv{(0,10)& (-6,10)} 
	?(.03)*\dir{-} ?(1)*\dir{-}; 
	"3";"2" **\crv{(6,-10)& (0,-10)} 
	?(.03)*\dir{-} ; 
	"1";"1E" **\dir{-}; 
	"3B";"3" **\dir{-};
	(-9,-5)*{a};
	(9,5)*{b};
	(2,0)*{i};
	\endxy
	\qquad = \qquad
	\xy 
	(0,10)*{}; (0,-10)*{} **\dir{-};
	(-2,10)*{b};
	(2,-10)*{a};
	\endxy
	\]
	\caption{$M_{ai}M^{ib} = \delta^{b}_{a}$}
	\label{fig:six}
\end{figure}

\begin{figure}
	\[
	\xy
	\vtwist~{(-5,5)}{(5,5)}{(-5,-10)}{(5,-10)}; 
	\vcross~{(-5,-10)}{(5,-10)}{(-5,-20)}{(5,-20)};
	(-5,7)*{a};
	(5,7)*{b};
	(-7,-10)*{i};
	(7,-10)*{j};
	(-5,-22)*{c};
	(5,-22)*{d};
	\endxy
	\qquad = \qquad
	\xy
	(-5,5)*{}; (-5,-20)*{} **\dir{-}; 
	(5,5)*{}; (5,-20)*{} **\dir{-}; 
	(-5,7)*{a};
	(5,7)*{b};
	(-5,-22)*{c};
	(5,-22)*{d};
	\endxy
	\]
	\caption{$ \overline{R^{ab}_{ij}} R^{ij}_{cd} = \delta^{a}_{c} \delta^{b}_{d}$}
	\label{fig:seven}
\end{figure}

\begin{figure}
	\[
	\xy 
	\vcross~{(-5,10)}{(5,10)}{(-5,0)}{(5,0)}; 
	(15,10)*{}; (15,0)*{} **\dir{-};
	\vcross~{(5,0)}{(15,0)}{(5,-10)}{(15,-10)}; 
	(-5,0)*{}; (-5,-10)*{} **\dir{-}; 
	\vcross~{(-5,-10)}{(5,-10)}{(-5,-20)}{(5,-20)}; 
	(15,-10)*{}; (15,-20)*{} **\dir{-}; 
	(-5,12)*{a};
	(5,12)*{b};
	(15,12)*{c};
	(-5,-22)*{e};
	(5,-22)*{f};
	(15,-22)*{g};
	(-7,-5)*{i};
	(8, 2)*{j};
	(8, -9)*{k};
	\endxy
	\qquad = \qquad
	\xy 
	\vcross~{(5,10)}{(15,10)}{(5,0)}{(15,0)}; 
	(-5,10)*{}; (-5,0)*{} **\dir{-};
	\vcross~{(-5,0)}{(5,0)}{(-5,-10)}{(5,-10)}; 
	(15,0)*{}; (15,-10)*{} **\dir{-}; 
	\vcross~{(5,-10)}{(15,-10)}{(5,-20)}{(15,-20)}; 
	(-5,-10)*{}; (-5,-20)*{} **\dir{-};
	(-5,12)*{a};
	(5,12)*{b};
	(15,12)*{c};
	(-5,-22)*{e};
	(5,-22)*{f};
	(15,-22)*{g};
	(8, 2)*{j};
	(8, -9)*{k};
	(17, -5)*{i};
	\endxy
	\]
	\caption{$R^{ab}_{ij} R^{jc}_{kg} R^{ik}_{ef} = R^{bc}_{ji} R^{aj}_{ek} R^{ki}_{fg}$}
	\label{fig:eight}
\end{figure}

\begin{figure}
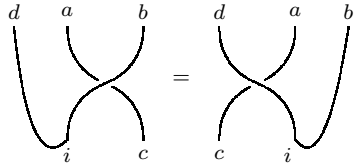

	\[
	\xy
	\vtwist~{(-5,5)}{(5,5)}{(-5,-10)}{(5,-10)}; 
	(-12,5)*{}="A"; 
	(-5,-10)*{}="B"; 
	"A"; "B" **\crv{(-10,-15)}; 
	\vcross~{(15,5)}{(25,5)}{(15,-10)}{(25,-10)}; 
	(32,5)*{}="C"; 
	(25,-10)*{}="D"; 
	"C"; "D" **\crv{(29,-15)}; 
	(10,-2)*{=};
	(-5,7)*{a};
	(5,7)*{b};
	(-5,-12)*{i};
	(5, -12)*{c};
	(-12, 7)*{d};
	(15,7)*{d};
	(25,7)*{a};
	(15,-12)*{c};
	(24, -12)*{i};
	(32, 7)*{b};
	\endxy 
	\]
	\caption{$M^{di} \overline{R}^{ab}_{ic} = R^{da}_{ci} M^{ib}$}
	\label{fig:nine}
\end{figure}

In the case of the Jones polynomial, we have all the algebra present to make the model. It is easiest to indicate the model for the bracket polynomial as given in \cite{state}: let cup and cap be given by the $2 \times 2$ matrix $M$, described above so that $M_{ij} = M^{ij}$. Let $R$ and $\overline{R}$ be given by the equations

$$
 R^{ab}_{cd} = A M^{ab} M_{cd} + A^{-1} \delta^{a}_{c} \delta^{b}_{d},
$$

$$
\overline{R^{ab}_{cd}} = A^{-1} M^{ab} M_{cd} + A \delta^{a}_{c} \delta^{b}_{d}.
$$

In general, the inverse of a matrix $R$ will be denoted by $\overline{R}$  throughout the discussion in the remainder of the paper.\\

The bracket is normalized so that the value of a circle is $-A^2 - A^{-2}$. In this specific case, we have the following matrix for $M$:

$$ M = \left[ \begin{array}{cc}
	0 & iA \\
-iA^{-1} & 0
\end{array} \right].
$$

This definition of the $R$ matrices exactly parallels the diagrammatic expansion of the bracket, and it is not hard to see, either by algebra or diagrams, that all the conditions of the model are met. Thus, this $R$ satisfies the Yang-Baxter equation. Other solutions to the Yang-Baxter equation give invariants distinct from the Jones polynomial.

\subsection{\bf Entanglement}

A unitary linear mapping $G : V \otimes V \rightarrow V \otimes V$ where $V$ is a two dimensional complex vector space and $G$ is some operator is said to be \textit{entangling} if there is a vector 

$$
\ket{\alpha \beta} = \ket{\alpha} \otimes \ket{\beta} \in V \otimes V
$$

\noindent such that $G \ket{\alpha \beta}$ is not decomposable as a tensor product of two qubits. Under these circumstances, one says that $G \ket{\alpha \beta}$ is \textit{entangled}. \\

\noindent \textbf{Example 2.1} A two-qubit  pure state 

$$
 \ket{\phi} = a\ket{00} + b \ket{01} + c \ket{10} + d \ket{11}
$$

\noindent is entangled exactly when $(ad-bc) \neq 0$ as proved in \cite{notsoshort}. It is easy to use this fact to check when a specific matrix is, or is not, entangling. 

\section{Unoriented State Models Given by Non-Entangling Operators}

\begin{figure}
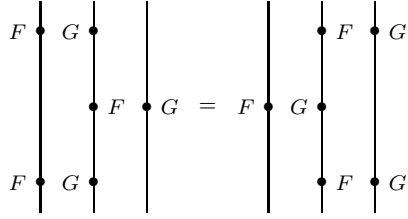

	\[
	\xy
	(0,0)*{}; (0,28)*{} **\dir{-}; 
	(7,0)*{}; (7,28)*{} **\dir{-}; 
	(-7,0)*{}; (-7,28)*{} **\dir{-};
	(-7,4)*+{\bullet}="4"; 
	(0,4)*+{\bullet}="6";
	(7,14)*+{\bullet}="4"; 
	(0,14)*+{\bullet}="6";
	(-7,24)*+{\bullet}="4"; 
	(0,24)*+{\bullet}="6";
	(-3,4)*{G};
	(-10,4)*{F};
	(3,14)*{F};
	(10,14)*{G};
	(-3,24)*{G};
	(-10,24)*{F};
	(30,0)*{}; (30,28)*{} **\dir{-}; 
	(37,0)*{}; (37,28)*{} **\dir{-}; 
	(23,0)*{}; (23,28)*{} **\dir{-};
	(37,4)*+{\bullet}="4"; 
	(30,4)*+{\bullet}="6";
	(23,14)*+{\bullet}="4"; 
	(30,14)*+{\bullet}="6";
	(37,24)*+{\bullet}="4"; 
	(30,24)*+{\bullet}="6";
	(33,4)*{F};
	(40,4)*{G};
	(27,14)*{G};
	(20,14)*{F};
	(33,24)*{F};
	(40,24)*{G};
	(15,14)*{=};
	\endxy
	\]
	\caption{This decomposition of the Yang-Baxter equation implies that $F^2 = xF, G^2 = tG.$}
	\label{fig:regdecomp}
\end{figure}

In \cite{amsshort}, the authors made use of the following theorem to characterize non-entangling operators. \\

\noindent \textbf{Theorem 3.1} \textit{Let $V$ be a finite-dimensional complex vector space, and $M \in GL(V \otimes V)$ be a non-entangling operator. Then there exist $A, B \in GL(V)$ such that either $M = A \otimes B$ or $M = (A \otimes B) \circ S$, where $S(x \otimes y) = y \otimes x$.}\\

The authors in \cite{amsshort} note that non-entangling operators are the invertible elements of $End(V \otimes V)$ which map product states to product states. The proof of this theorem is given in \cite{amsshort}. We call the two cases of this theorem the \textit{product case} for $M= A \otimes B$ and the \textit{swap case} for $M=(A\otimes B) \circ S$. In the following, we discuss state summation models for link invariants with respect to the two cases. 

\subsection{\bf The Product Case}

\begin{figure}
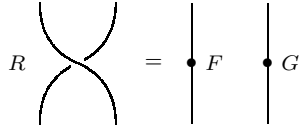

	\[
	\xy
	\vtwist~{(5,16)}{(-5,16)}{(5,0)}{(-5,0)}; 
	(15,16)*{}; (15,0)*{} **\dir{-}; 
	(25,16)*{}; (25,0)*{} **\dir{-}; 
	(15,8)*+{\bullet}="4";
	(25,8)*+{\bullet}="6";
	(18,8)*{F};
	(28,8)*{G};
	(10,8)*{=};
	(-8,8)*{R};
	\endxy
	\]
	\caption{Topological relations for the product case. $\overline{R}$ similarly decomposes to $\overline{F}$ and $\overline{G}$ on the identity.}
	\label{fig:regrelate}
\end{figure}

We now examine state summation models constructed given that $R = F \otimes G$ as shown in Figure~\ref{fig:regrelate}. The goal is to show that {\it when we decompose the $R$ matrix in this fashion the resulting state summation leads to a trivial invariant.} In order to accomplish this aim, we assume that $R$ has the form given above, and analyse the effect that this must have on the cup and cap evaluations. This means that we do not 
actually write cup and cap matrices in doing the analysis. We deduce the form of the invariant from the given conditions, and show that it must be a trivial invariant. Thus we go back to the basic diagrammatic 
restrictions that are imposed by Figure~\ref{fig:six}, Figure~\ref{fig:seven}, Figure~\ref{fig:eight}, Figure~\ref{fig:nine} and deduce conditons that are needed to produce an invariant. This same method of analysis is used 
throughout the rest of the paper.\\

Our methods are based on the state summation models for knots and links described in \cite{state}. In the arguments given below, we assume that a state summation model is given, using this $R$-matrix, and we deduce enough 
aspects of its structure to conclude that it is a trivial invariant.\\

From the Yang-Baxter equation as shown in Figure~\ref{fig:regdecomp}, we can deduce the fact that $F^2 = x F$ and $G^2 = tG$. As $F$ and $G$ are invertible, then $F = x I$ and $G = t I$, where $I$ is the identity. Therefore, $ R = s I $ where $s=xt$. This fact is also demonstrated in \cite{amsshort}. We now conclude that $R = s I$ and $\overline{R} = \overline{s} I$, where $\overline{s} = s^{-1}$. The relations are 

$$\langle \raisebox{-0.3\height}{\includegraphics[width=0.6cm, height=0.5cm]{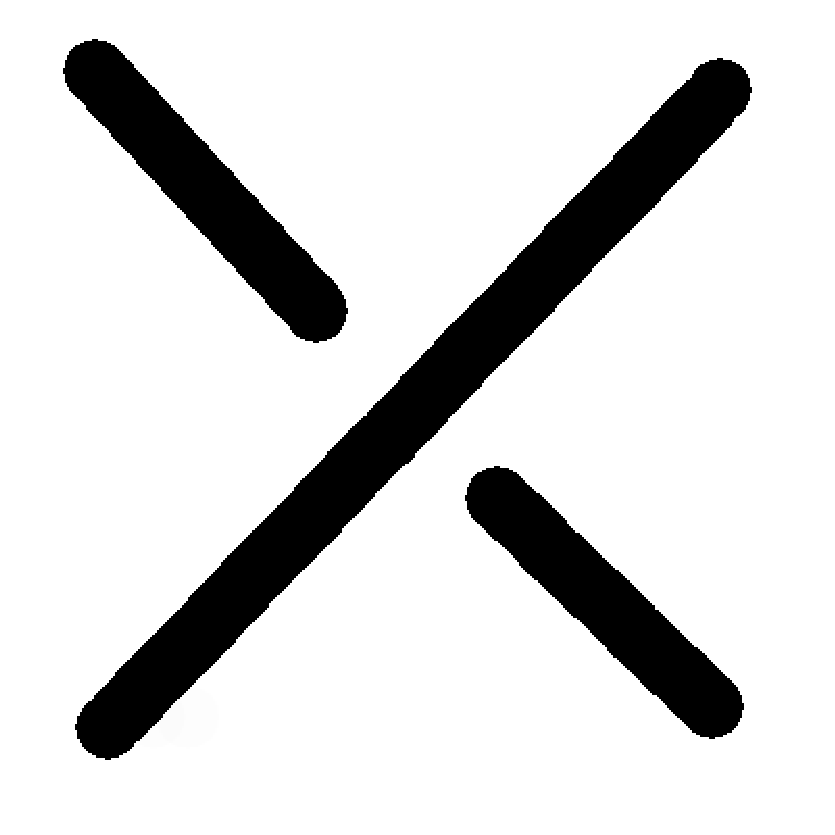}} \rangle = \overline{s} \langle \raisebox{-0.25\height}{\includegraphics[width=0.45cm, height=0.5cm]{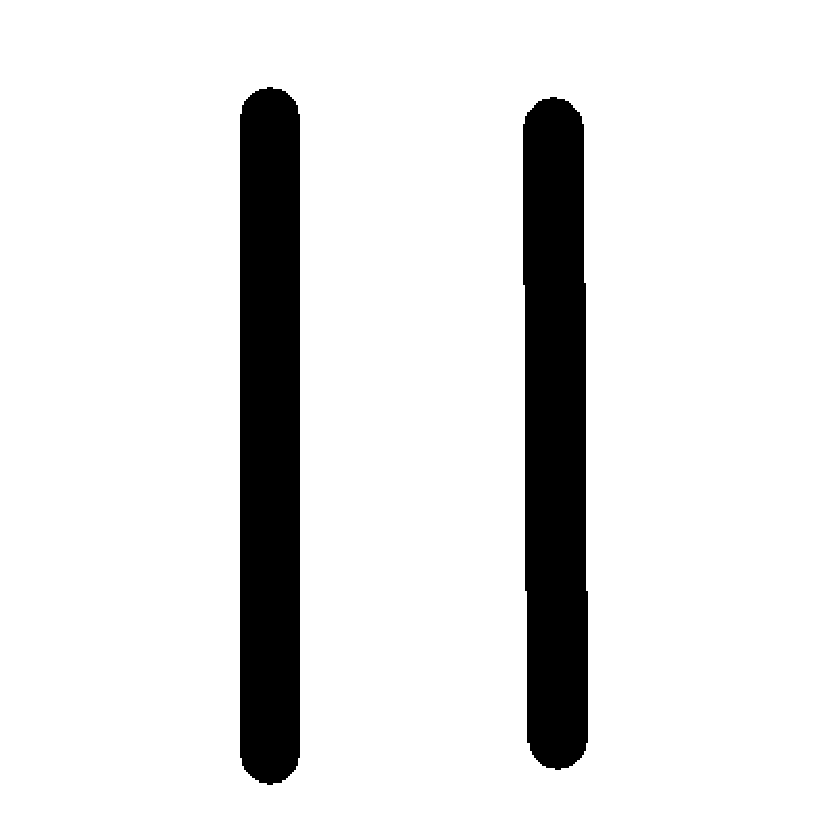}} \rangle$$

$$\langle \raisebox{-0.3\height}{\includegraphics[width=0.6cm, height=0.5cm]{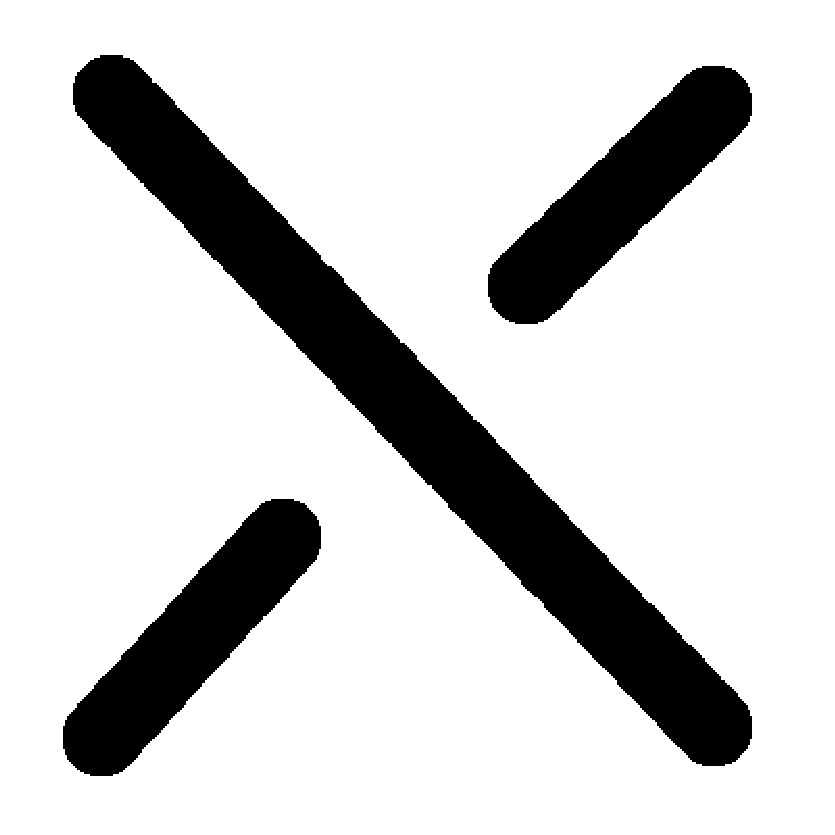}} \rangle = s \langle \raisebox{-0.25\height}{\includegraphics[width=0.45cm, height=0.5cm]{identity.eps}} \rangle$$

\noindent We use the following lemmas to construct an invariant from the state summation given by the above relations. \\

\noindent \textbf{Lemma 3.3}
$$\langle \raisebox{-0.25\height}{\includegraphics[width=0.45cm, height=0.5cm]{identity.eps}} \rangle =
	s^2 \langle  \raisebox{-0.2\height}{\includegraphics[width=0.35cm, height = 0.35cm]{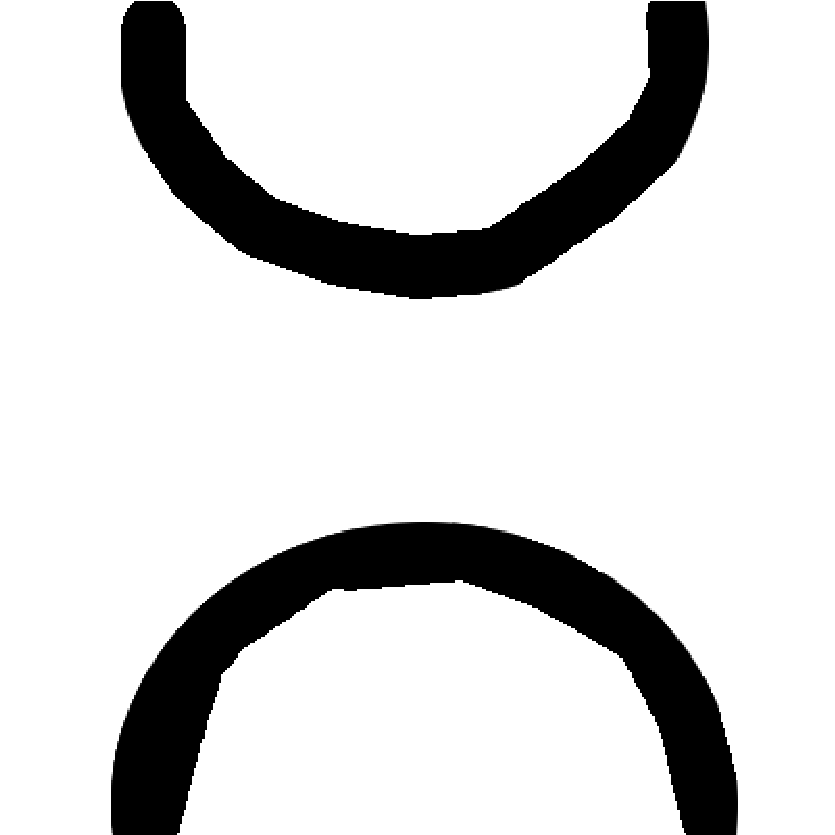}} \rangle. $$
\begin{proof}
	Note that the relation
	$$\langle \raisebox{-0.3\height}{\includegraphics[width=0.6cm, height=0.5cm]{rcross.eps}} \rangle =
	 \langle  \raisebox{-0.3\height}{\includegraphics[width=0.6cm, height=0.5cm]{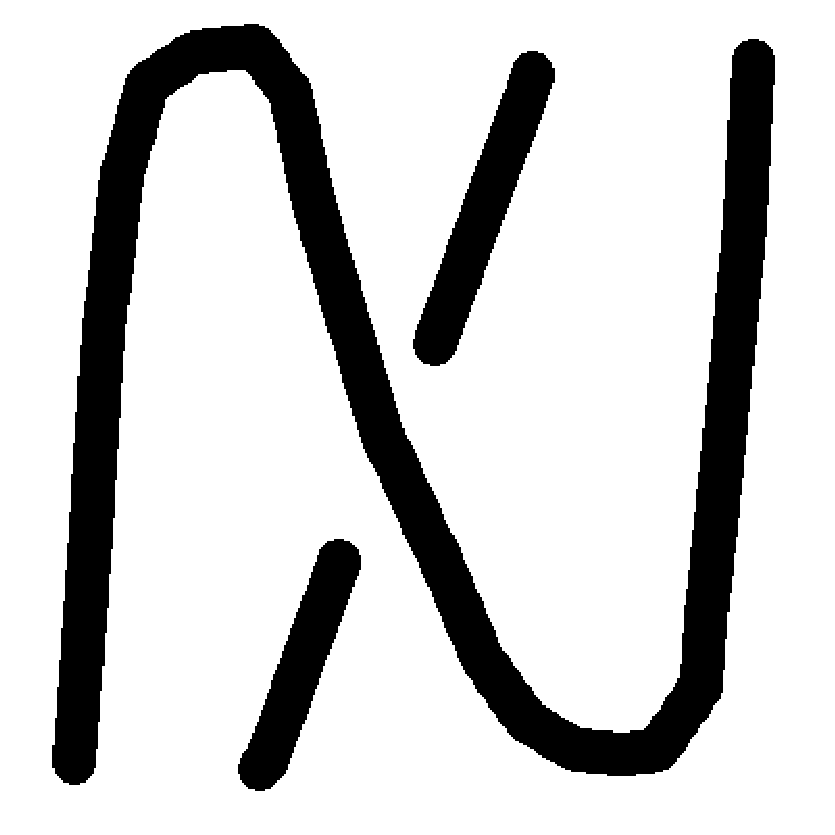}} \rangle$$
is independent of the particular choice of cup or cap matrices. This is analogous to twisting $\overline{R}$. By applying the smoothings associated to $R$ and $\overline{R}$, we arrive at the following:
	$$\overline{s} \langle \raisebox{-0.2\height}{\includegraphics[width=0.45cm, height=0.5cm]{identity.eps}} \rangle =
	s \langle  \raisebox{-0.2\height}{\includegraphics[width=0.65cm, height = 0.45cm]{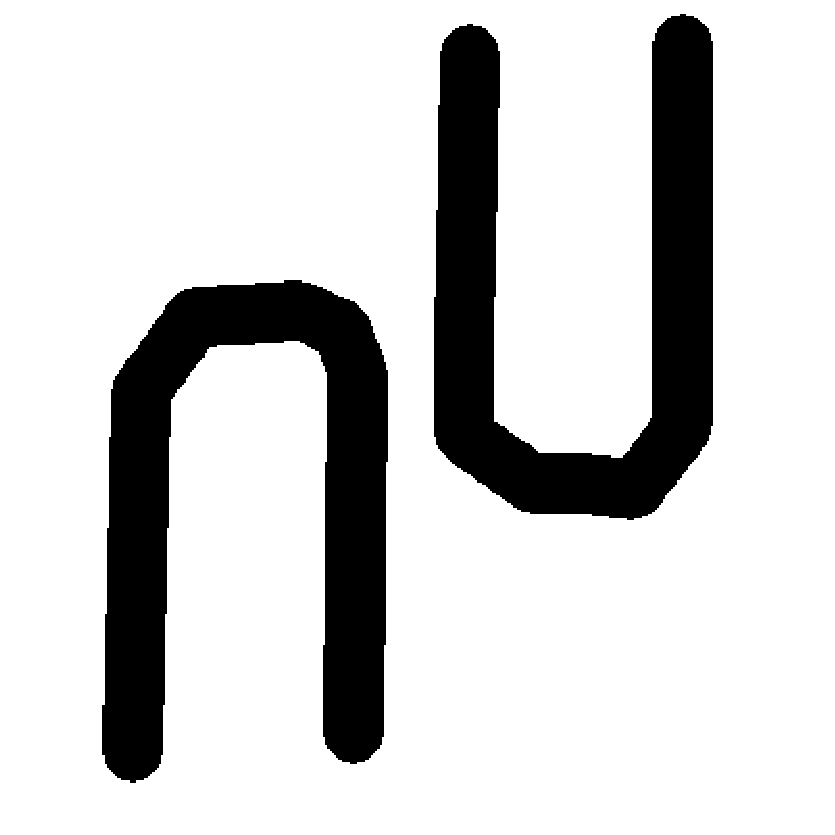}} \rangle,$$
	$$\overline{s} \langle \raisebox{-0.2\height}{\includegraphics[width=0.45cm, height=0.5cm]{identity.eps}} \rangle =
	s \langle  \raisebox{-0.2\height}{\includegraphics[width=0.38cm, height = 0.35cm]{capcap.eps}} \rangle,$$
	$$\langle \raisebox{-0.25\height}{\includegraphics[width=0.45cm, height=0.5cm]{identity.eps}} \rangle =
	s^2 \langle  \raisebox{-0.2\height}{\includegraphics[width=0.38cm, height = 0.35cm]{capcap.eps}} \rangle. $$ \qed
\end{proof}	

\noindent \textbf{Corollary 3.3.1}
	$$ \langle \raisebox{-0.25\height}{\includegraphics[width=0.5cm, height=0.5cm]{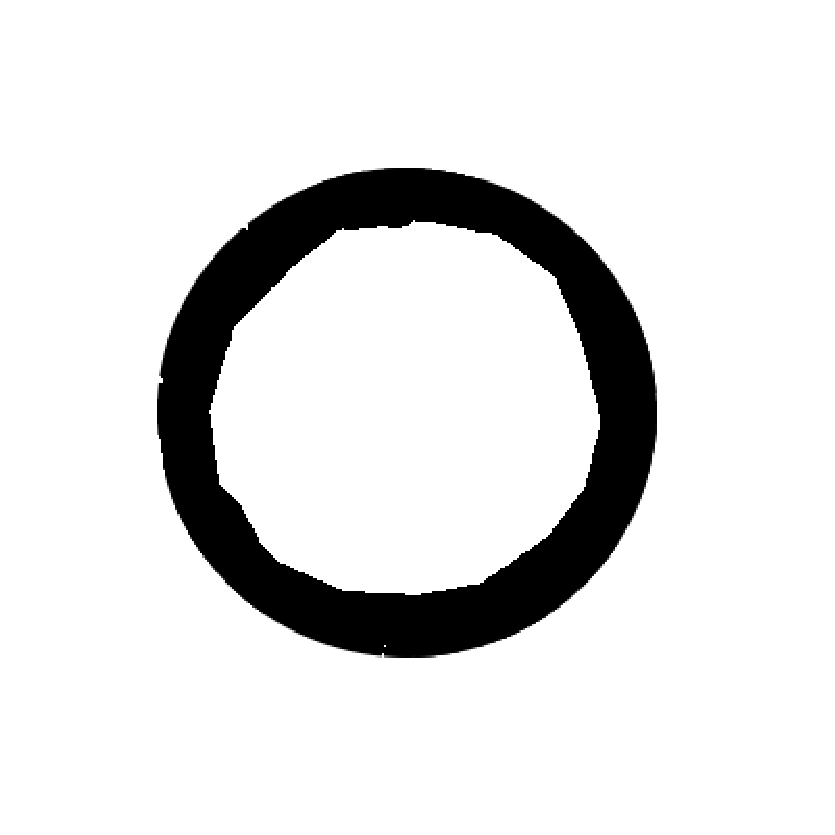}} \rangle =
	s^2 \langle  \raisebox{-0.2\height}{\includegraphics[width=0.55cm, height = 0.45cm]{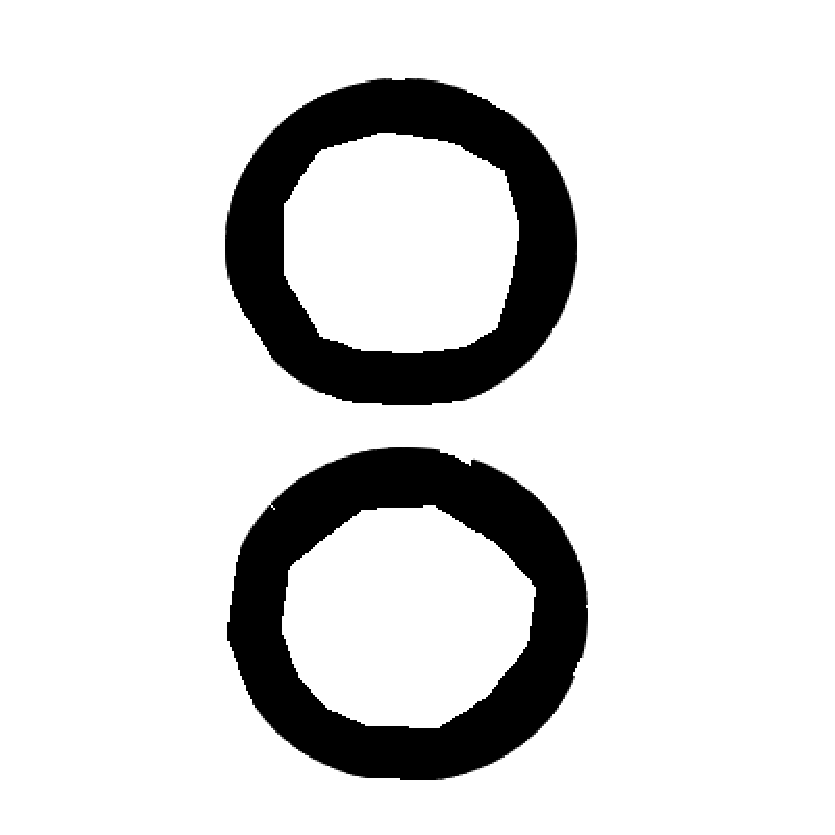}} \rangle. $$

\noindent \textbf{Corollary 3.3.2}
	$$ \langle \raisebox{-0.25\height}{\includegraphics[width=0.5cm, height=0.5cm]{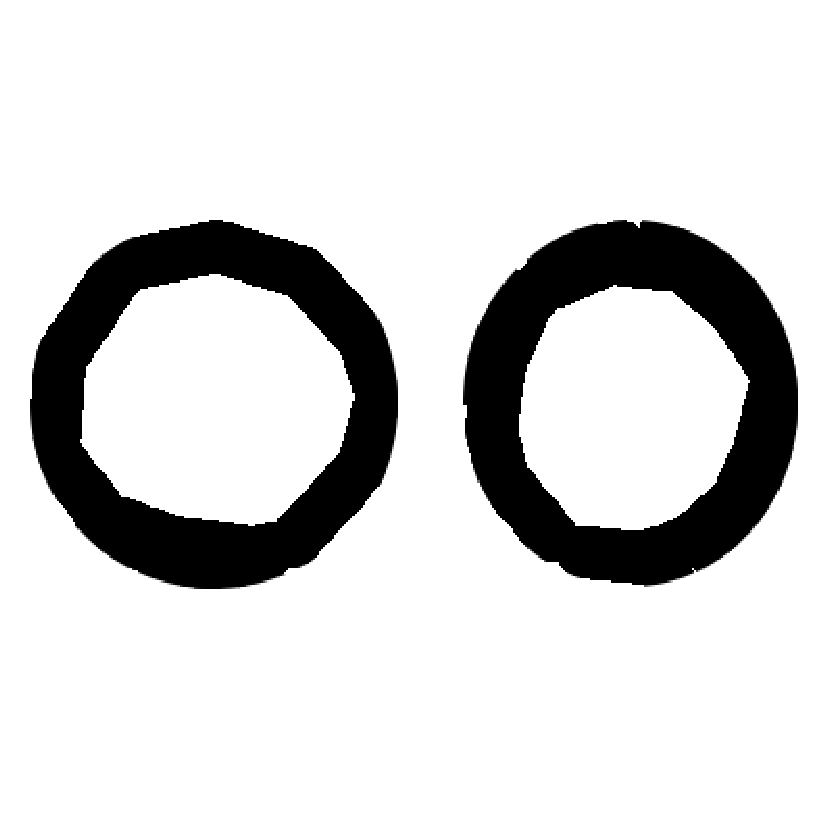}} \rangle =
	s^2 \langle  \raisebox{-0.2\height}{\includegraphics[width=0.55cm, height = 0.45cm]{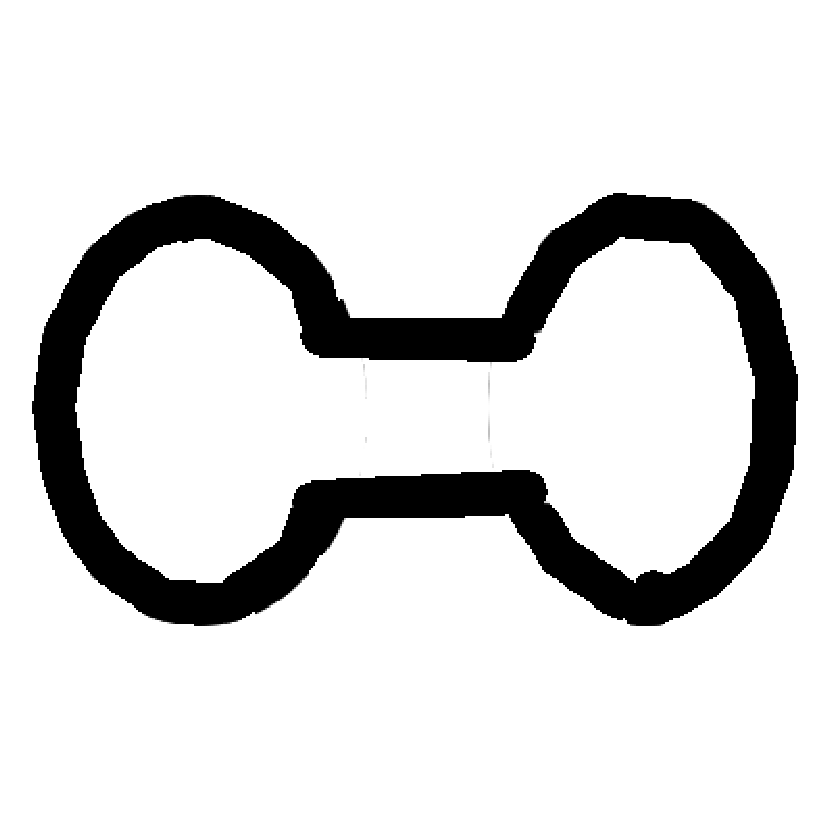}} \rangle.$$

\noindent Setting the value of the circle equal to $\delta$, we have that $\delta s^2 = 1$ and $s^2 = \delta$. We now arrive at the fact that $s^4 = 1$.

\noindent \textbf{Lemma 3.4} \textit{(The Second Reidemeister Move)} Invariance of the state summation under the second Reidemeister move follows from the formal properties we have given so far.
\begin{proof}
By applying our smoothing to the following diagram and then using Lemma 3.3 we get
	$$ \langle \raisebox{-0.25\height}{\includegraphics[width=0.5cm, height=0.5cm]{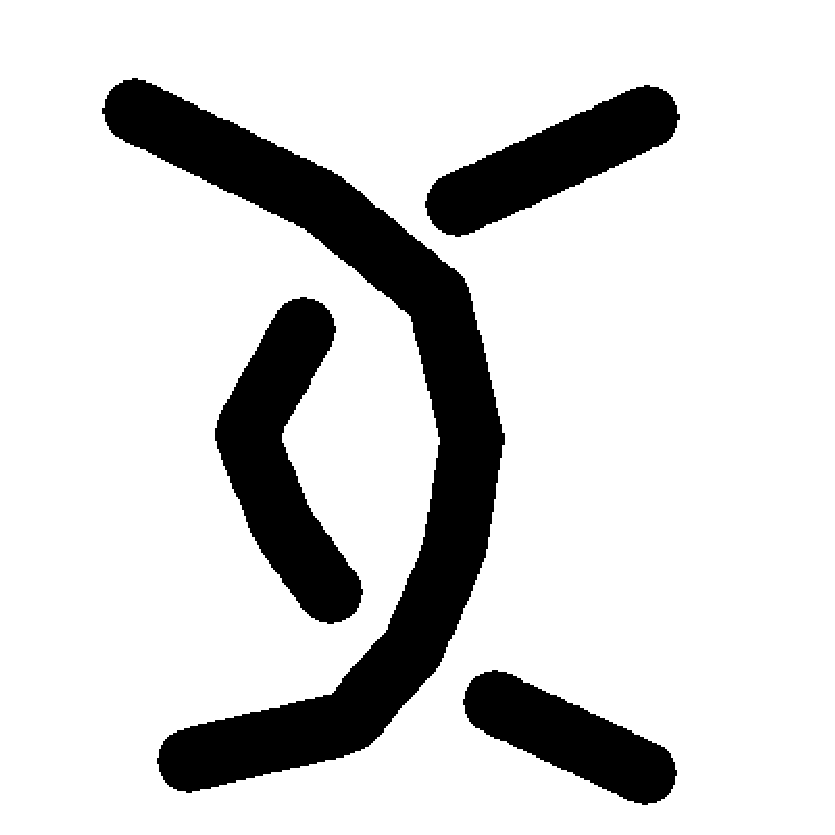}} \rangle =
	s^2 \langle  \raisebox{-0.25\height}{\includegraphics[width=0.5cm, height = 0.4cm]{capcap.eps}} \rangle =  \langle \raisebox{-0.25\height}{\includegraphics[width=0.5cm, height=0.5cm]{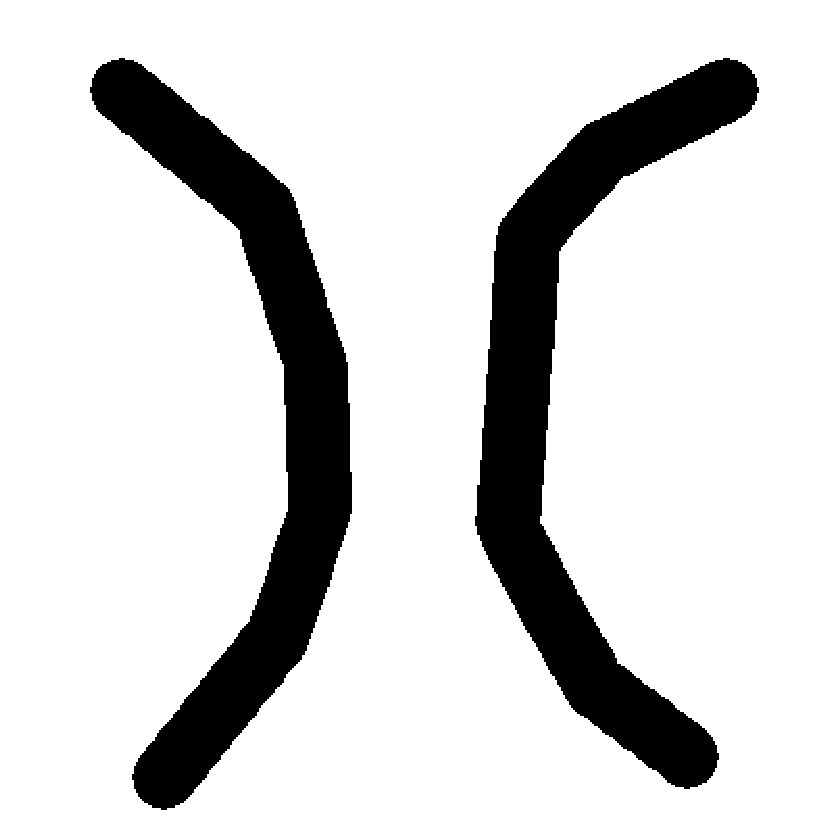}} \rangle.$$ \qed
\end{proof}

\noindent \textbf{Lemma 3.5} \textit{(The Third Reidemeister Move)} Invariance of the state summation under the third Reidemeister move follows from the formal properties we have given so far.
\begin{proof}
The third Reidemeister move immediately follows by replacing one crossing by a smoothing as shown just before Lemma 3.3,  and then using Lemma 3.4. \qed
\end{proof}
	
\noindent \textbf{Lemma 3.6} \textit{(The First Reidemeister Move)} The state sum multiplies by $s$ for positive curls and by $\bar{s}$ for negative curls.

\begin{proof}
Since the relations are
$$\langle \raisebox{-0.3\height}{\includegraphics[width=0.6cm, height=0.5cm]{rcross.eps}} \rangle = \overline{s} \langle \raisebox{-0.25\height}{\includegraphics[width=0.45cm, height=0.5cm]{identity.eps}} \rangle,$$
$$\langle \raisebox{-0.3\height}{\includegraphics[width=0.6cm, height=0.5cm]{cross.eps}} \rangle = s \langle \raisebox{-0.25\height}{\includegraphics[width=0.45cm, height=0.5cm]{identity.eps}} \rangle,$$
we can apply them to the curls.

%$$ \langle \raisebox{-0.25\height}{\includegraphics[width=0.4cm, height=0.4cm]{loop.eps}} \rangle =
%	s \langle  \raisebox{-0.2\height}{\includegraphics[width=0.4cm, height = 0.35cm]{line.eps}} \rangle. $$1
%	$\overline{R}$ will then have a writhe of $1$ in the loop below,
%	$$ \langle \raisebox{-0.25\height}{\includegraphics[width=0.4cm, height=0.4cm]{iloop.eps}} \rangle =
%	s \langle  \raisebox{-0.2\height}{\includegraphics[width=0.4cm, height = 0.35cm]{circline.eps}} \rangle, $$
	
	$$ \langle \raisebox{-0.25\height}{\includegraphics[width=0.4cm, height=0.4cm]{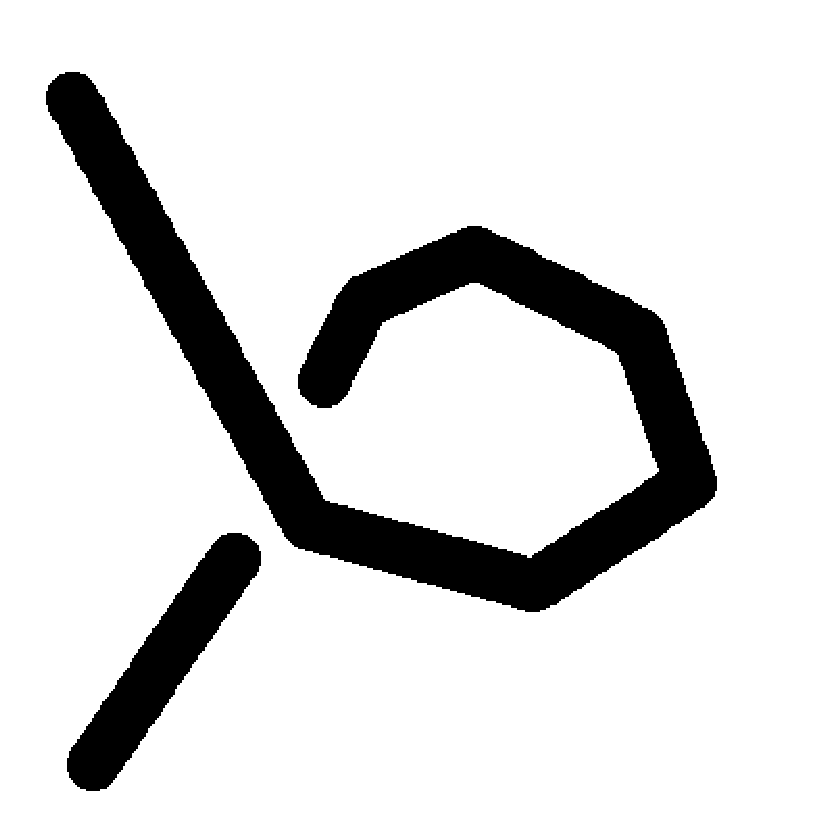}} \rangle = s \langle \raisebox{-0.25\height}{\includegraphics[width=0.4cm, height=0.4cm]{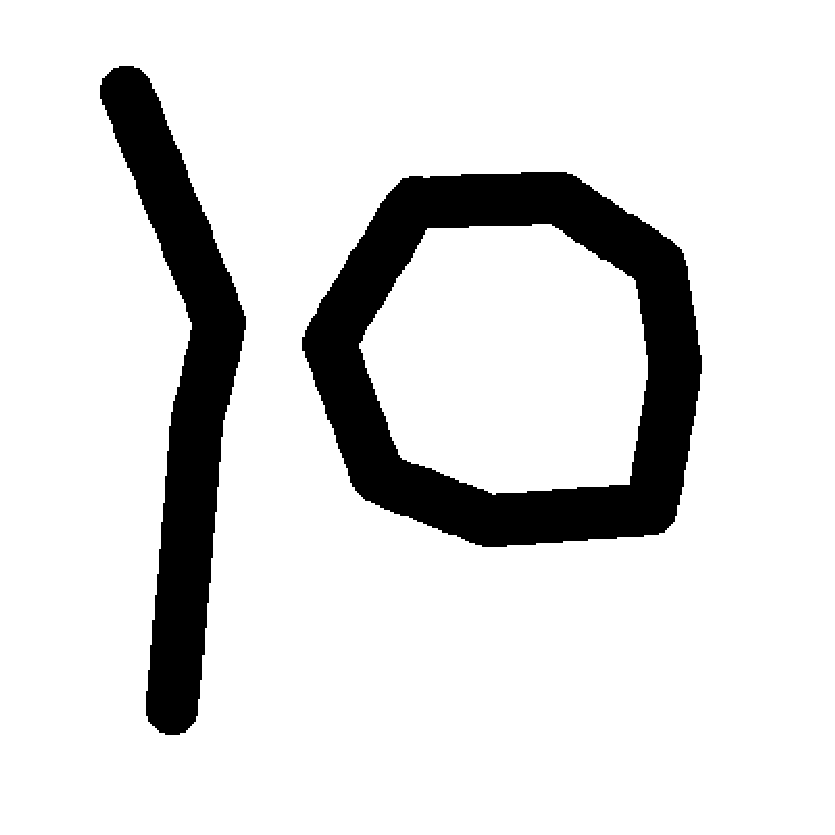}} \rangle
	= s \delta \langle  \raisebox{-0.2\height}{\includegraphics[width=0.4cm, height = 0.35cm]{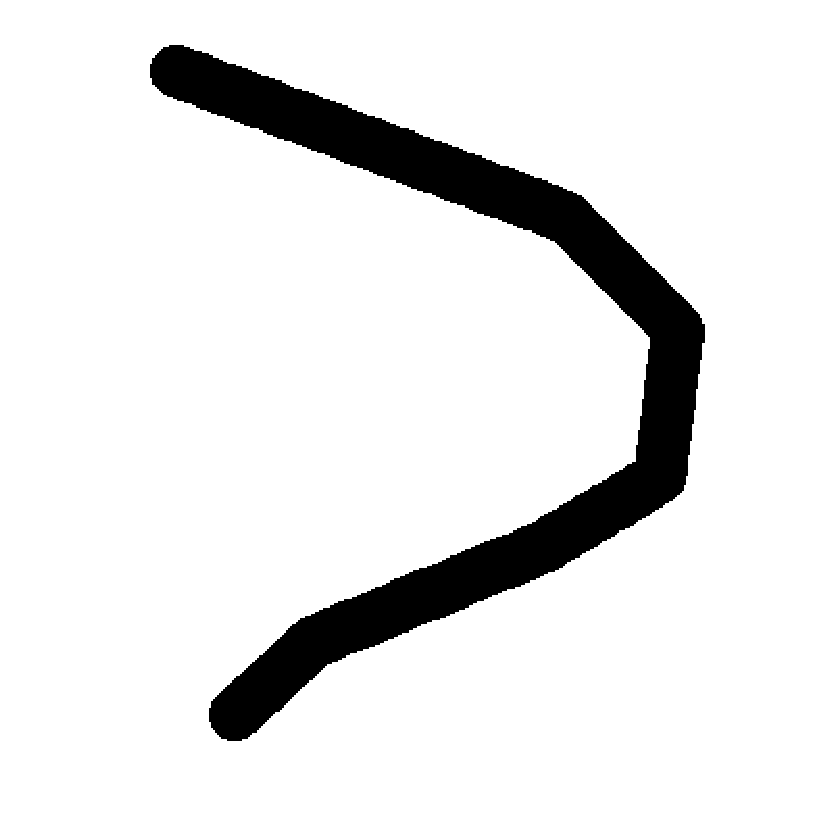}} \rangle 
	= s^3 \langle  \raisebox{-0.2\height}{\includegraphics[width=0.4cm, height = 0.35cm]{line.eps}} \rangle =  \overline{s} \langle  \raisebox{-0.2\height}{\includegraphics[width=0.4cm, height = 0.35cm]{line.eps}} \rangle.$$ 
	The other relation follows in the same fashion.  
	\qed
\end{proof}

\noindent \textbf{Theorem 3.2} \textit{The quantum state summation given by  $R = F \otimes G$ is a trivial invariant of unoriented knots.}

\begin{proof} 

In order to get an ambient isotopy invariant $f_K$ for knots, we would need to compensate for the extra factors that arise from performing the first Reidemeister move. We accomplish this via writhe-normalization as in \cite{state}. For a knot $K$ we define $f_K$ by the equation

$$f_K = s^{-w(K)} \langle K \rangle. $$

In order to use this formula, orient the knot diagram and then smooth it in an oriented way at every crossings. The result of this smoothing is the collection of {\it Seifert Circles} for the diagram. Let $SC(K)$ denote the number of Seifert circles in $K.$ Using the results above including the writhe compensation it is easy to see that each crossing contributes $s^{2sgn(c)}$ where $sgn(c)$ denotes the sign of the crossing. The factor of $2$ 
occurs because both sign of crossing and smoothing of crossing each contribute $s^{sgn(c)}.$ From this it follows that $$f_K = s^{-2 wr(K)} \delta^{SC(K)} = \delta^{-wr(K) + SC(K)}.$$
The Lemma in the Appendix to this paper shows that $$ - wr(K) + SC(K)  \equiv 1\pmod 2.$$
Therefore, since $\delta^2 =1,$ we conclude that  $f_{K} = \delta$ for all knots $K.$
This completes the proof of the Theorem. 
\qed

\end{proof}

\begin{figure}
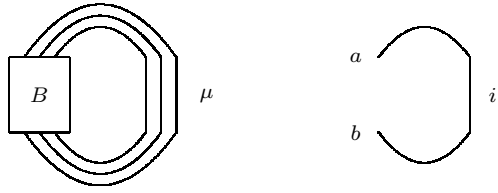

	\[ 
	\xy
	(5,0)*{}; (5,10)*{} **\dir{-};
	(-3,0)*{}; (-3,10)*{} **\dir{-};
	(5,10)*{}; (-3,10)*{} **\dir{-};
	(5,0)*{}; (-3, 0)*{} **\dir{-};
	(17,10)*{}; (17,0)*{} **\dir{-};
	(19,10)*{}; (19,0)*{} **\dir{-};
	(15,10)*{}; (15,0)*{} **\dir{-};
	(3,10)*{}="A"; 
	(15,10)*{}="B"; 
	"A"; "B" **\crv{(9,18)}; 
	(1,10)*{}="C"; 
	(17,10)*{}="D"; 
	"C"; "D" **\crv{(9,21)};
	(-1,10)*{}="E"; 
	(19,10)*{}="F"; 
	"E"; "F" **\crv{(9,24)};
	(3,0)*{}="G"; 
	(15,0)*{}="H"; 
	"G"; "H" **\crv{(9,-8)}; 
	(1,0)*{}="I"; 
	(17,0)*{}="J"; 
	"I"; "J" **\crv{(9,-11)};
	(-1,0)*{}="K"; 
	(19,0)*{}="L"; 
	"K"; "L" **\crv{(9,-14)};
	(1,5)*{B};
	(23,5)*{\mu};
	\endxy
	\qquad \qquad \qquad
	\xy
	(15,10)*{}; (15,0)*{} **\dir{-};
	(3,10)*{}="A"; 
	(15,10)*{}="B"; 
	"A"; "B" **\crv{(9,18)}; 
	(3,0)*{}="G"; 
	(15,0)*{}="H"; 
	"G"; "H" **\crv{(9,-8)}; 
	(0,10)*{a};
	(0,0)*{b};
	(18,5)*{i};
	\endxy
	\]
	\caption{In braid closures the enhancement operator $\mu$ must correspond to a cup and a cap.}
	\label{fig:mu}
\end{figure}

\subsection{\bf The Swap Case}
In Figure~\ref{fig:swaprelate} we show the form of the braiding operator for the unoriented swap case. We begin this section by analyzing the state sum models for operators of this form.
For an unoriented knot or link diagram $K$ in Morse form, we will let the invariant of regular isotopy associated with this braiding operator be denoted by $Invar(K).$ \\

\noindent \textbf{Theorem} \textit {The state sum model $Invar(K)$ for $K$  in the unoriented swap case produces only trivial invariants for knots (links of one component).}

\begin{proof}

First note that via Figure~\ref{yangbaxterswap} we have that the Yang-Baxter equation for $R$ implies 
that $FG = GF$ in the swap case where $F$ and $G$ appear in $R$ as in Figure~\ref{fig:swaprelate}. Then from Figure~\ref{swing} we conclude that $\overline{F} = F$ and $\overline{G} = G$ so that 
$F^2 = 1 = G^2.$ The Figure~\ref{swing} shows that we can slide $F$ and $G$ over maxima and minima in the diagram leaving them unchanged. This means that in a knot diagram we can collect all algebra on a diagram as a single product along a given arc.  Since the number of $F$'s equals the number of crossings, and the number of $G$'s equals the number of crossings, we have that the algebraic expression can be written in the form
$F^n G^n$ where $n$ is the number of crossings in the knot diagram. Since $F^2 = G^2 = 1,$ we can take the exponent $n$ modulo two. We also know that $S \equiv n + 1$ ( mod $2$) where $S$ denotes the number of Seifert circuits in the knot diagram. This follows from the Lemma in the Appendix to this paper. Furthermore, the Whitney degree of the underlying plane curve of the diagram is congruent modulo two to $S$ \cite{FKT}.
It follows that $F^n G^n = FG$ when the Whitney degree is even and  $F^n G^n = I$ when  the Whitney degree is odd. Taking into account the fact that every Reidemeister type one move contributes $1$ to the Whitney degree and contributes  $FG$ to the algebraic part of the evaluation, we see that the evaluation of any knot diagram is the same as the evaluation of a corresponding unknot diagram with the same writhe and Whitney degree. This completes the proof 
that the knot invariant cannot distinguish any knot from the unknot. \qed 

\end{proof}

It remains to discuss the possibility that the invariant could detect a link. For the purpose of this discussion
we shall take the cup and cap operators for this model to be identity operators. That is, we shall assume that $M_{ab} = \delta_{ab}$ and $M^{ab} = \delta^{ab}$ where $\delta_{ab} = \delta^{ab} = 1$ if and only if
$a=b$ and $\delta_{ab} = \delta^{ab} = 0$ when $a \ne b.$ 
From Figure~\ref{hopf} we see that if the evaluation of a loop of Whitney degree one, labeled with an algebraic expression $\alpha,$ is denoted by
$Tr(\alpha)$, then the evaluation of the Hopf Link as shown in this figure is $Tr(FG)Tr(FG) = Tr(FG)^2.$ The corresponding evaluation of an unlink is $Tr(I)^2.$ We will now give an explicit example for $F$ and $G$ where these two evaluations differ, showing that an invariant in the unoriented swap case can detect linking even though the Yang-Baxter operator is not entangling. Consider the matrices $F$ and $G$ shown below.\\

$$F= \left[ \begin{array}{ccc}
0 & 0 & 1 \\
0 & 1 & 0 \\
1 & 0 & 0 
\end{array} \right],
G= \left[ \begin{array}{ccc}
1 & 0 & 0 \\
0 & -1 & 0 \\
0 & 0 & 1 
\end{array} \right],
FG= \left[ \begin{array}{ccc}
0 & 0 & 1 \\
0 & -1 & 0 \\
1 & 0 & 0
\end{array} \right].
$$\\

It is easy to verify that $F^2 = G^2 = I$ and that $FG = GF.$ The state sum model will use $Tr(\alpha) = Trace(\alpha)$ where {\it Trace} denotes standard matrix trace.  This gives a consistent state model. Note that both
$F$ and $G$ are symmetric matrices and that this corresponds to the invariance of the slide over maxima and minima in Figure~\ref{swing}. We then have (since these are $3 \times 3$ matrices) that 
$Tr(I) = 3$ while $Tr(FG) = -1.$ Thus $(Tr(FG))^2 = 1$ while $(Tr(I))^2 = 9$ and so this invariant detects the Hopf Link.\\

\noindent {\bf Remark.} Note that the result of doing a first Reidemeister move for the invariant under discussion is to multiply the algebra element on the component on which the move occurs by $FG.$ See Figure~\ref{curl}.
Since the algebra on a given component is either $FG$ or the identity $I,$ we see that the result of a first Reidemeister move is to switch the value of the invariant on this component from $Tr(FG) = -1$ to $Tr(I) = 3$ or from 
$3$ to $-1.$ The simplest way to use the invariant as an invariant of ambient isotopy is to use the fact: {\it Two links with the same Whitney degree and writhe (for each component)  are regularly isotopic if and only if they are ambient isotopic.} See \cite{OnKnots}. In this way we can prepare diagrams for comparison. This is how we know that the Hopf Link as shown in Figure~\ref{hopf} is shown to be non-trivial by this invariant. The two components of the Hopf Link diagram used in the calculation give results identical to two disjoint circles for the unlink.\\

We have the following result:\\

\noindent \textbf{Theorem} \textit {The state sum model $Invar(K)$ for links $K$ of two components can detect the modulo two linking number  of any link of two components and is non-trivial for 
links of odd linking number and trivial for links of even linking number.}\\ 

\begin{proof}
The proof follows from the discussion above and an easy analysis of the products of algebra elements that occur on the link components. \qed
\end{proof}

\noindent{\bf Remark.} We underline the fact that we have constructed a state sum invariant of knots and links, based on a {\it non-entangling} Yang-Baxter operater ( of swap type) that can detect the Hopf Link. 
This shows that the state sum models in this swap case have a similar relationship with linking and quantum entanglement as do the enhanced Yang-Baxter operators using Markov trace as in Section 3 of 
\cite{amsshort}, where an example of the detection of the Hopf link is given in a different way. The state sum that we have described here does not fit into the braiding form with enhanced Yang-Baxter operator that is 
used in \cite{amsshort}, but our state sum is indeed based on a Yang-Baxter operator. The examples in both cases show that non-entangling Yang-Baxter operators can detect non-trival topological linking.\\

\begin{figure}
	\[ 
	\xy 
	\vcross~{(-5,0)}{(5,0)}{(-5,-20)}{(5,-20)};
	(20,0)*{}; (30,-20)*{} **\dir{-};
	(30,0)*{}; (20,-20)*{} **\dir{-};
	(0,-22)*{R};
	(27.4,-15)*+{\bullet};
	(22.6,-15)*+{\bullet};
	(20,-15)*{F};
	(30,-15)*{G};
	(12.5,-10)*{=};
	\vtwist~{(50,0)}{(60,0)}{(50,-20)}{(60,-20)};
	(75,0)*{}; (85,-20)*{} **\dir{-};
	(85,0)*{}; (75,-20)*{} **\dir{-};
	(82.4,-5)*+{\bullet};
	(77.6,-5)*+{\bullet};
	(74,-5)*{\overline{F}};
	(86,-5)*{\overline{G}};
	(67.5,-10)*{=};
	(55,-22)*{\overline{R}};
	\endxy
	\]
	\caption{Topological relations for the swap case.}
	\label{fig:swaprelate}
\end{figure}

\begin{figure}
     \begin{center}
     \begin{tabular}{c}
     \includegraphics[width=6cm]{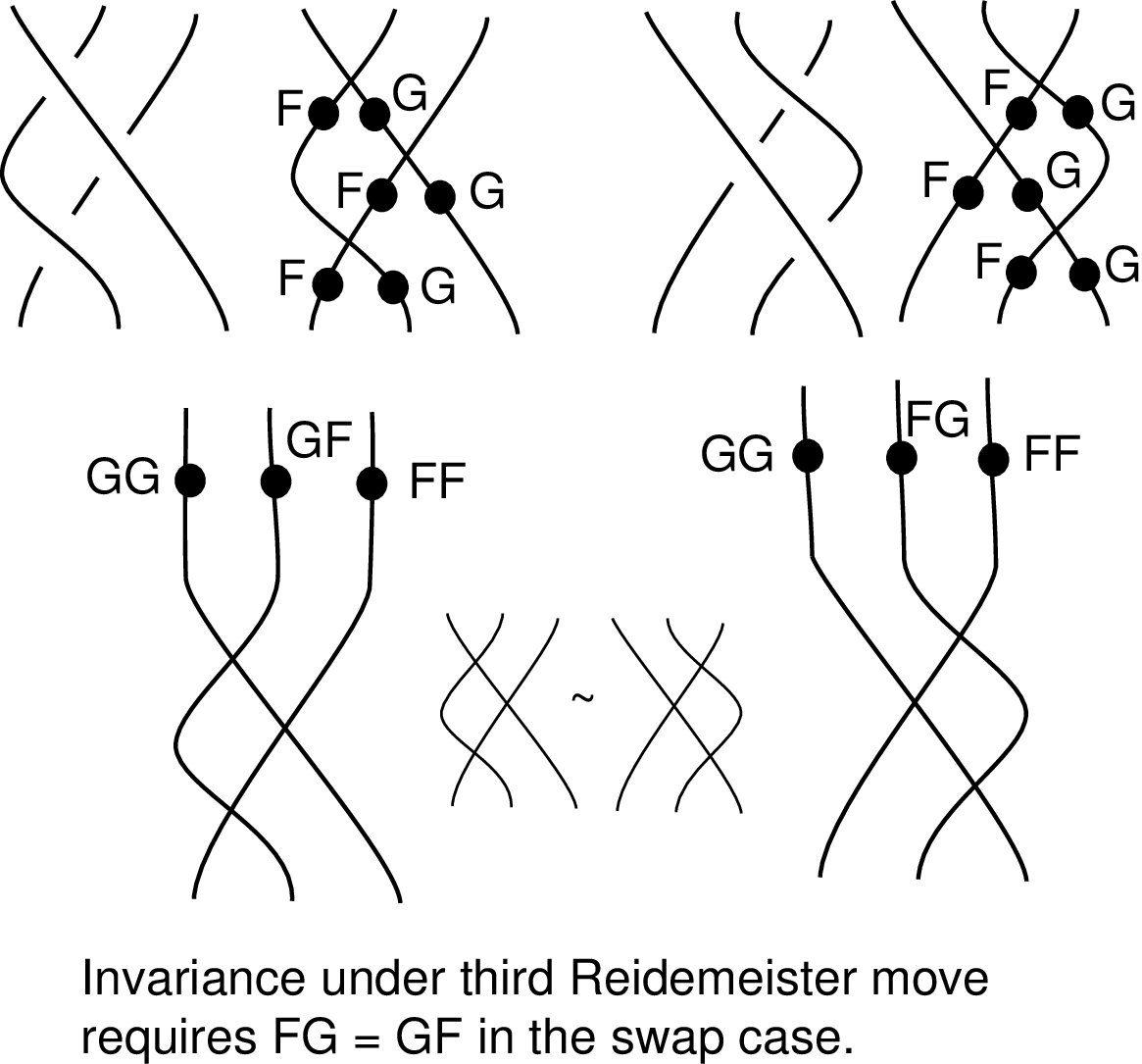}
     \end{tabular}
     \caption{\bf Third Reidemeister Move Implies $FG = GF.$}
     \label{yangbaxterswap}
\end{center}
\end{figure}

\begin{figure}
     \begin{center}
     \begin{tabular}{c}
     \includegraphics[width=6cm]{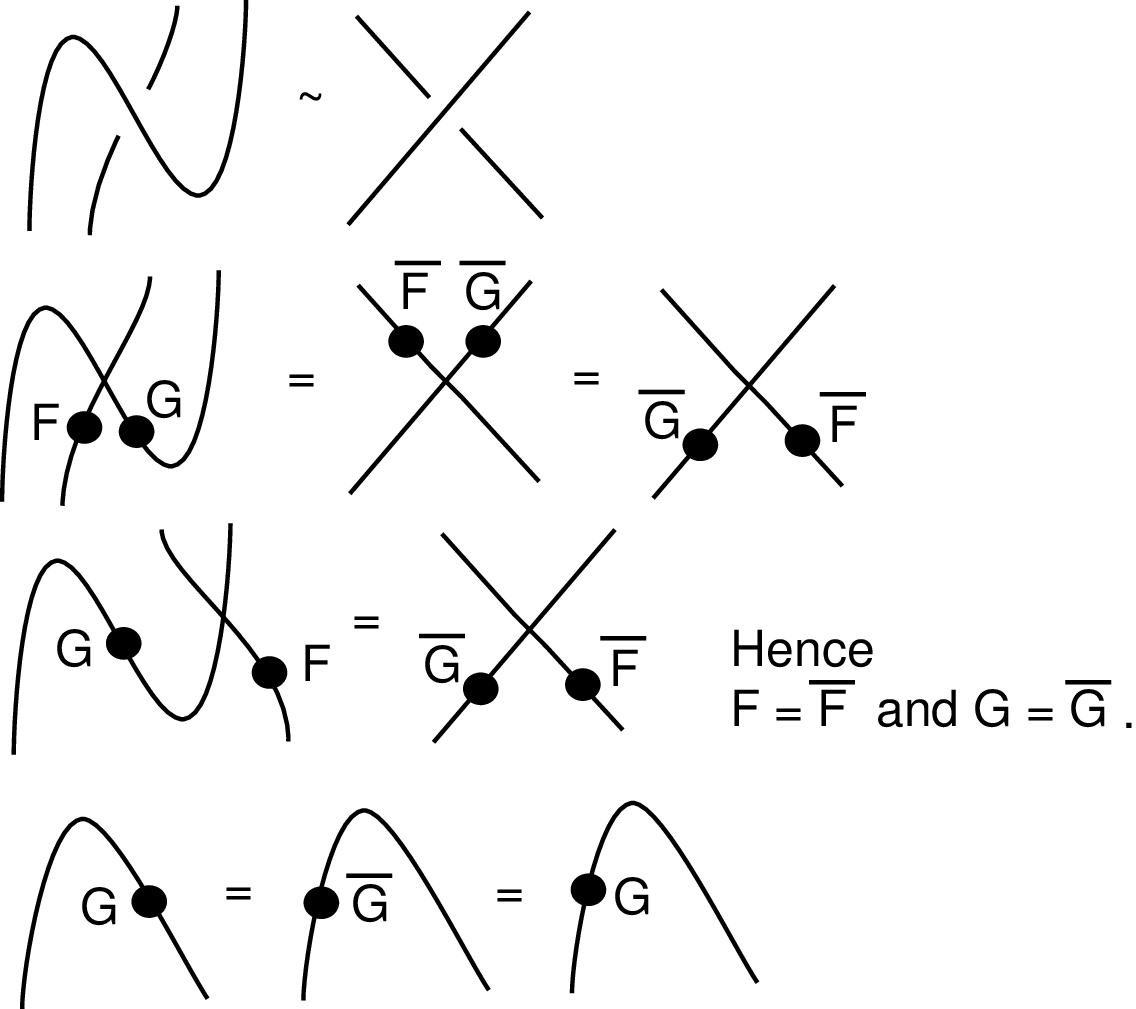}
     \end{tabular}
     \caption{\bf Swing Move Implies $FF = GG = I.$}
     \label{swing}
\end{center}
\end{figure}

\begin{figure}
     \begin{center}
     \begin{tabular}{c}
     \includegraphics[width=6cm]{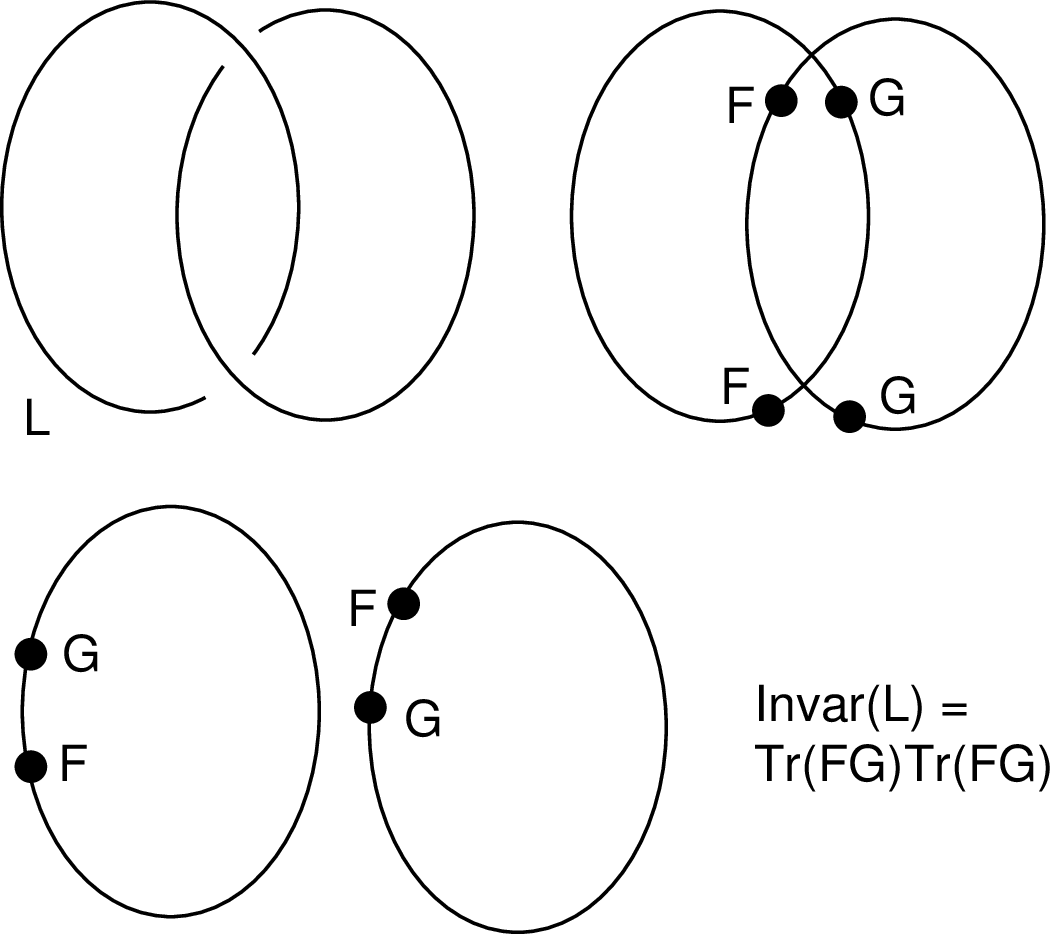}
     \end{tabular}
     \caption{\bf Invariant of the Hopf Link.}
     \label{hopf}
\end{center}
\end{figure}

\begin{figure}
     \begin{center}
     \begin{tabular}{c}
     \includegraphics[width=6cm]{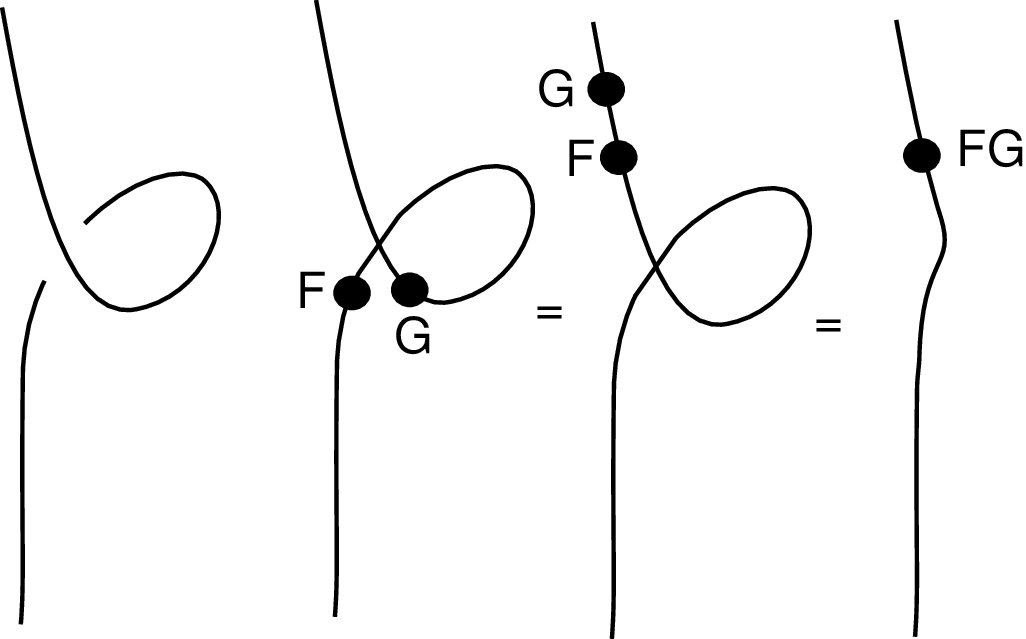}
     \end{tabular}
     \caption{\bf Behaviour under a Curl.}
     \label{curl}
\end{center}
\end{figure}

We now compare the methods of  \cite{amsshort} and our methods in this swap case. We rely heavily on quantum link invariants in this part of the paper. The decomposition of an $R$ matrix in the the swap $S \circ (F \otimes G)$ is represented topologically for both $R$ and $\overline{R}$ in Figure~\ref{fig:swaprelate}. \\

The proof in \cite{amsshort} relies on the use of enhanced Yang-Baxter operators defined below. \\

\noindent \textbf{Definition 3.7} \textit{Let $V$ be a finite-dimensional complex Hilbert space, $R \in GL(V \otimes V)$ a Yang-Baxter operator, and $\mu \in End(V)$. If $R$ commutes with $\mu \otimes \mu$ and
	$$ Tr_{2} (R \cdot \mu \otimes \mu) = Tr_{2} (\overline{R} \cdot \mu \otimes \mu) = \mu$$
then we say that the pair $\mathcal{R} = (R, \mu)$ is an enhanced Yang-Baxter operator. In that case, given any braid $b$ we define
	$$I_{\mathcal{R}} (b) = Tr[\rho_{n}^{R} (b) \cdot \mu^{\otimes n}].$$}

For a braid closure, the enhancement operator is analogous to a product of cup and cap matrices $M^{bi}$ and $M_{ai}$  taking $b$ to be a point on the strand at the top of the braid, $a$ to be the corresponding point on the strand at the bottom of the braid, and $i$ to be the point in the middle dividing the cup and the cap. This analogy is shown in Figure~\ref{fig:mu}. Note also that $M^{bi} = M_{bi}$ and $M^{ab}$ is the inverse of $M_{ab}$. Therefore, note that $\mu = M^{bi}M_{ai} = ( M M^{\top})_{ab}$. Given $M = (M^{ab})$ is a cup and $N = (M_{ab})$ is a cap (they are inverses), in general we have that 

$$ M = \left[ \begin{array}{cc}
a & b \\
c & d
\end{array} \right],
\qquad
M^{-1} = N = \frac{1}{\Delta} \left[ \begin{array}{cc}
d & -b \\
-c & a
\end{array} \right],
$$
	
\noindent with $\Delta = ad - bc$ as $M$ and $N$ must be inverses of one another. Thus

$$ \mu = (MN^{\top})_{ab} = \frac{1}{\Delta} \left[ \begin{array}{cc}
a & b \\
c & d
\end{array} \right]
\left[ \begin{array}{cc}
d & -c \\
-b & a
\end{array} \right]
= \frac{1}{\Delta} \left[ \begin{array}{cc}
ad-b^2 & -ac+ab \\
cd-bd & -c^2 + ad
\end{array} \right].
$$

In \cite{amsshort}, the authors found the following invertible $\mu$ for constructing an invariant. We will show that this $\mu$ cannot be obtained by the cup and cap construction.
(It is also the case that in \cite{amsshort} non-invertible $\mu$ are considered and shown to be unnecessary, but that is not the issue here.)

$$ \mu = \left[ \begin{array}{cc}
1 & 0 \\
0 & -1
\end{array} \right].
$$

\noindent If we try to set our derivation of $\mu$ from cup and cap operators equal to theirs, we then get matrix solutions that have a determinant of zero.
\begin{center}
$$\begin{array}{rcl} ad - b^{2} & = & \Delta, \\ -ac + ab & = & 0, \\ cd-bd & = & 0, \\ c^2 - ad & = & \Delta. \end{array}$$
\end{center}
Solving the above system of equations implies that $d = \frac{c^2}{a}$, which makes the determinant of the matrix zero. Therefore, in general we cannot extend the models used in \cite{amsshort} to quantum link invariants. This shows that we cannot always construct an analogous state summation model for the unoriented swap case. Nevertheless, we have shown that all properties relative to non-entangling Yang-Baxter operators
happen in the same way in both of these categories of invariants.\\

\noindent \textbf{Remark 3.8} \textit{We believe this is the first time an explicit difference has been shown between the Markov trace and state summation methods of constructing invariants.}

\section{Oriented State Models Given by Non-Entangling Operators}

We now express the above arguments in the oriented case for quantum invariant state summations. The results of \cite{amsshort} generalize easily in this case. 
See \cite{stat} for an account of oriented state sum models for link invariants, based on solutions to the Yang-Baxter equation. These models are essentially the same as the unoriented models, but have orientations associated with the crossings, cups, and caps. We begin with the simple product decomposition of $R$.

\subsection{\bf The Product Case}

\textbf{Theorem 4.1} \textit{The state summation model given by $R = F \otimes G$ is trivial for oriented knots.}

\begin{proof}

For the oriented product case, we begin with the following equations, as given in Section 3: 

$$\langle \raisebox{-0.2\height}{\includegraphics[width=0.5cm, height=0.4cm]{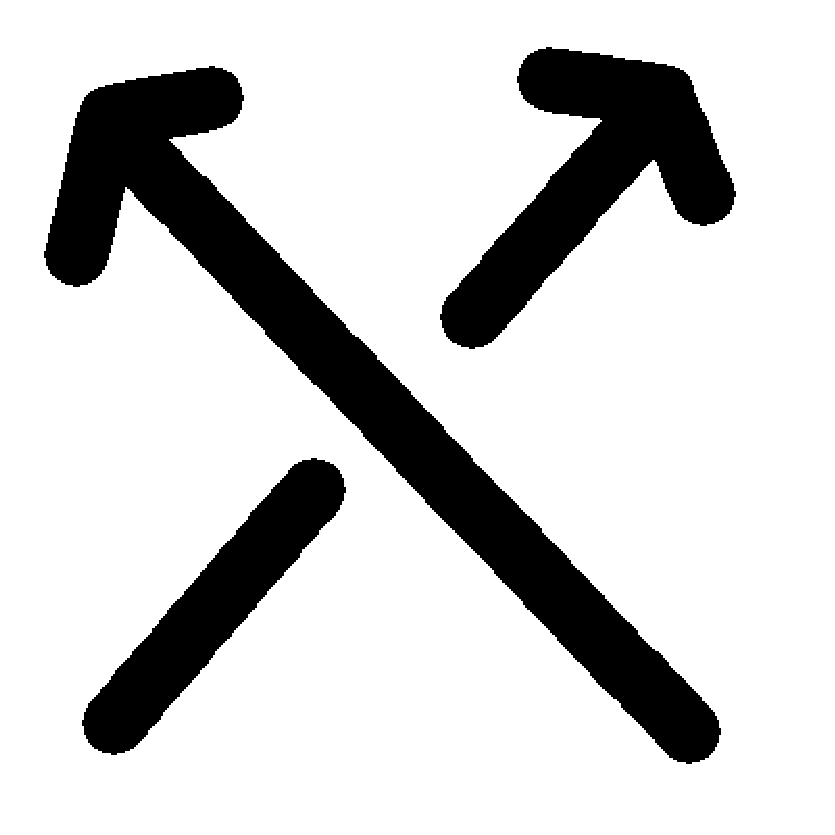}} \rangle =
s \langle \raisebox{-0.2\height}{\includegraphics[width=0.5cm, height=0.5cm]{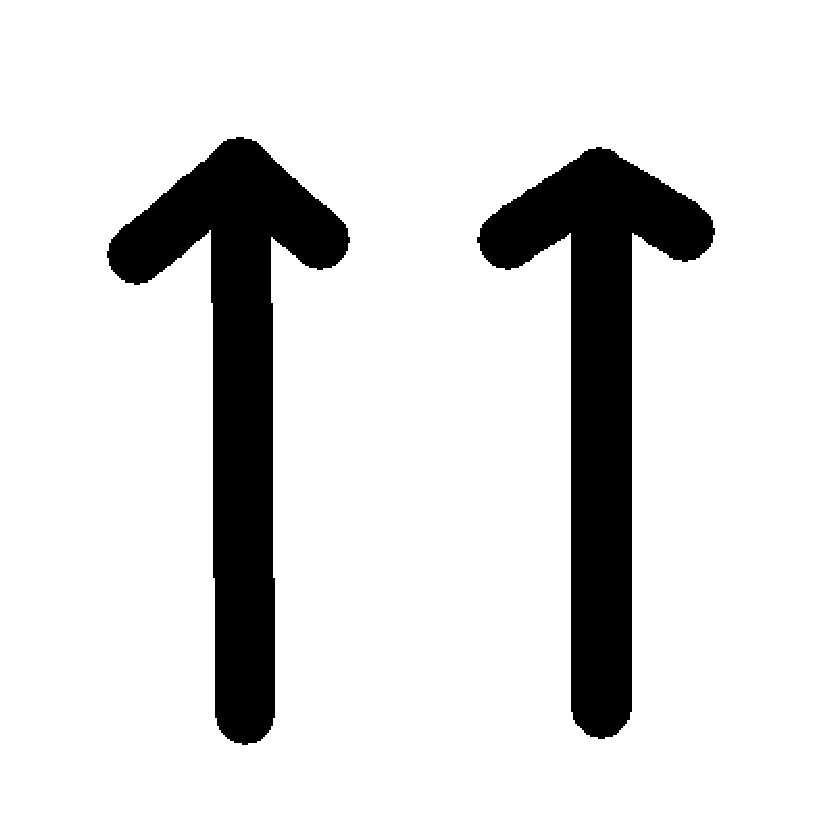}} \rangle,$$

$$\langle \raisebox{-0.2\height}{\includegraphics[width=0.5cm, height=0.4cm]{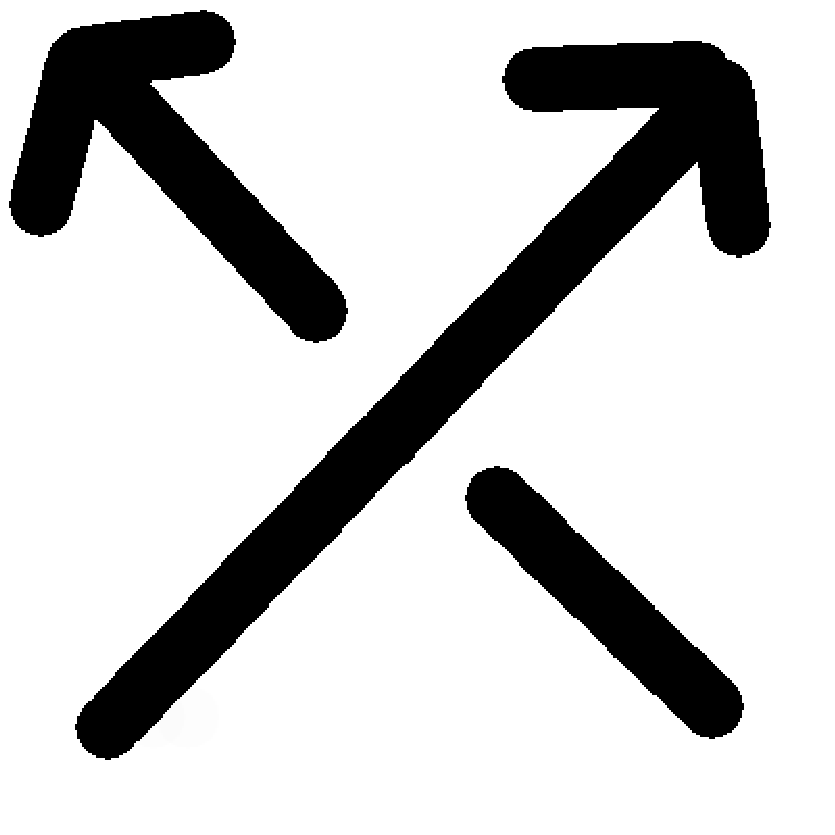}} \rangle =
\overline{s} \langle \raisebox{-0.2\height}{\includegraphics[width=0.5cm, height=0.5cm]{identityorient.eps}} \rangle.$$

We find the value of the loop by examining how the model behaves on the oriented curls. We apply our oriented smoothing relations to determine what factors arise from each curl:

$$\langle \raisebox{-0.2\height}{\includegraphics[width=0.45cm, height=0.45cm]{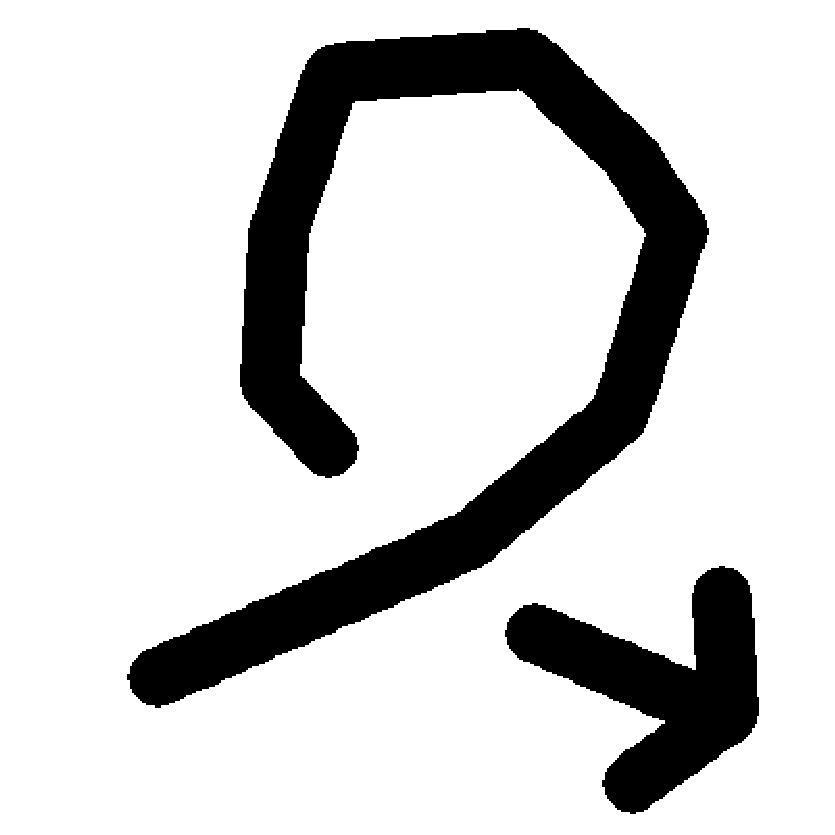}} \rangle = \overline{s} \langle \raisebox{-0.2\height}{\includegraphics[width=0.45cm, height=0.45cm]{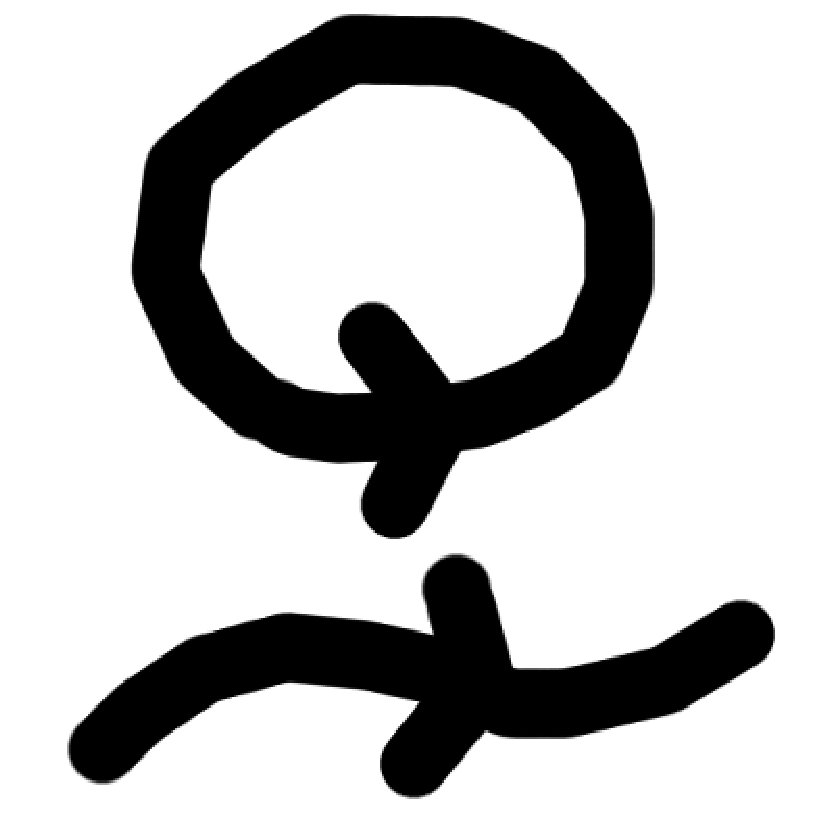}} \rangle = \overline{s} \delta \langle \raisebox{-0.2\height}{\includegraphics[width=0.45cm, height=0.45cm]{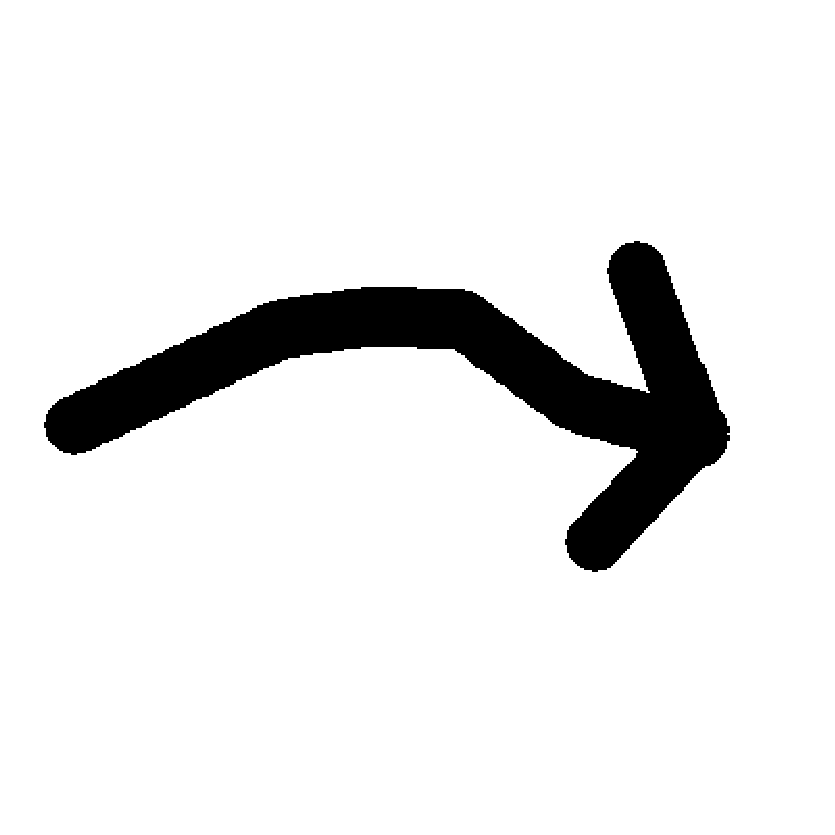}} \rangle$$

$$\langle \raisebox{-0.2\height}{\includegraphics[width=0.45cm, height=0.45cm]{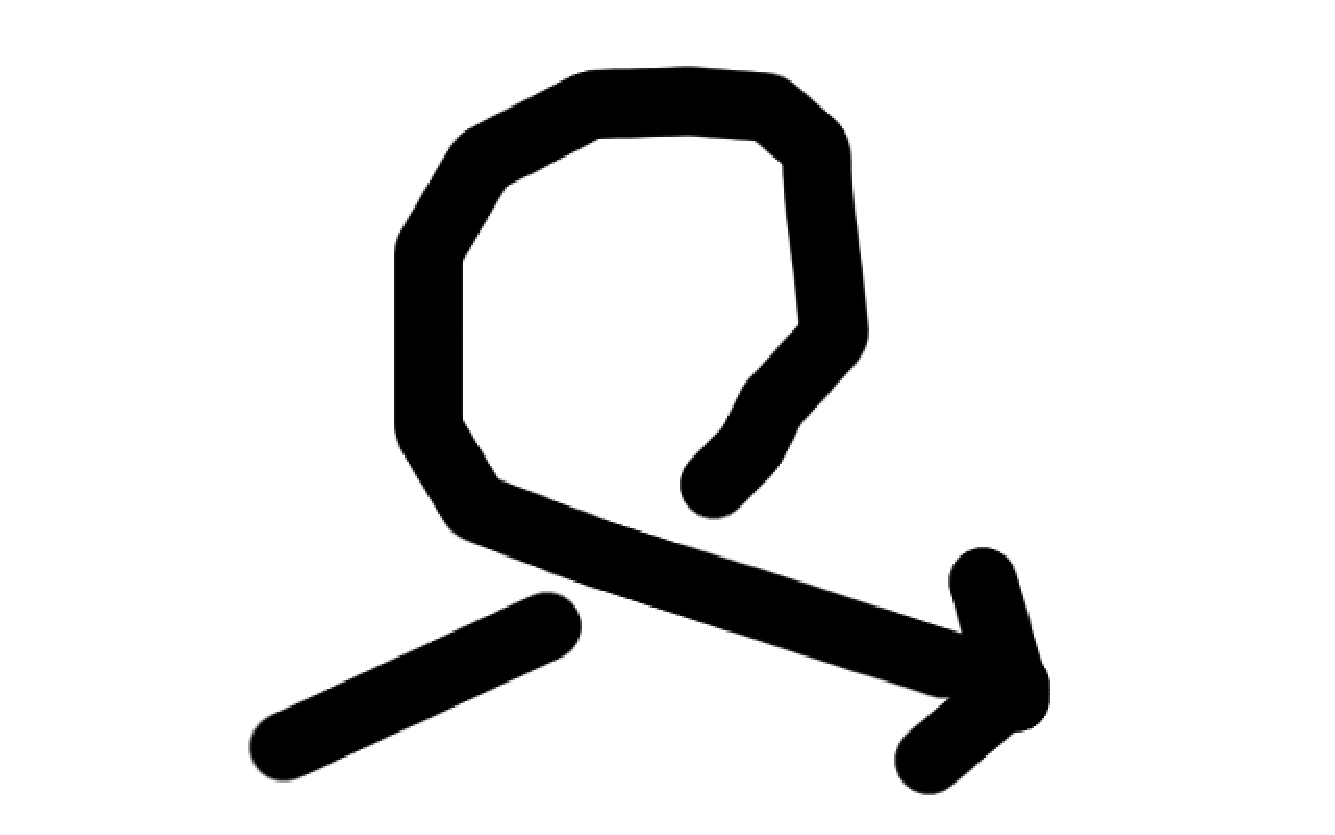}} \rangle = s \langle \raisebox{-0.2\height}{\includegraphics[width=0.45cm, height=0.45cm]{orientloopline.eps}} \rangle = s \delta \langle \raisebox{-0.2\height}{\includegraphics[width=0.45cm, height=0.45cm]{lineorient.eps}} \rangle,$$

\noindent taking $\delta$ to be the value of the oriented loop. In order to satisfy the relation $\langle \raisebox{-0.2\height}{\includegraphics[width=0.80cm, height=0.45cm]{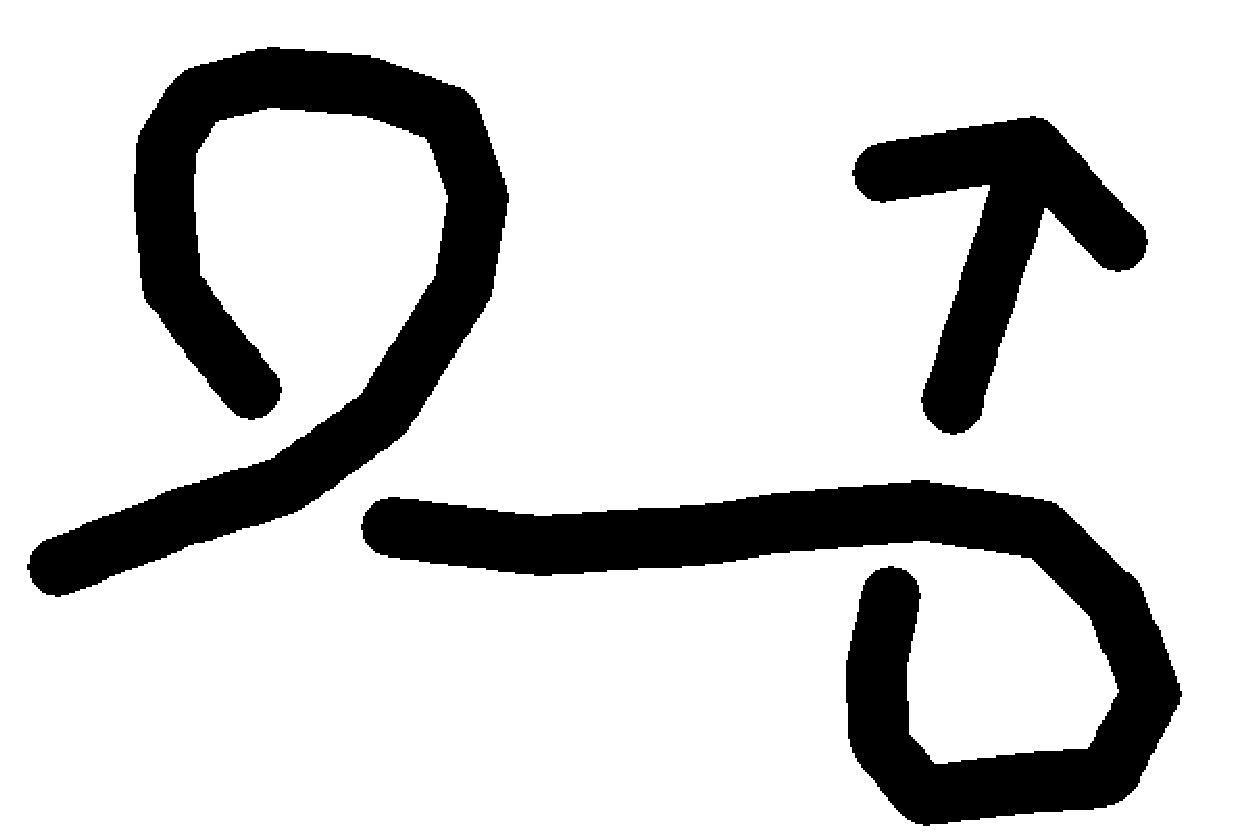}} \rangle = \langle \raisebox{-0.2\height}{\includegraphics[width=0.45cm, height=0.45cm]{lineorient.eps}} \rangle$ we need that $(\overline{s}\delta)(s\delta)= \delta^2 = 1$. Thus, the absolute value of an oriented loop is equal to one. To construct a model that is invariant under Reidemeister I we must multiply by a writhe corrective factor, so our final invariant has the form

$$f_k = (s \delta)^{-w(K)} \langle K \rangle,$$

\noindent as $s \delta$ is the value that comes from using our oriented smoothing relations on a positive writhe looped strand. \\

Note that if $K$ is a given diagram, then we have an explicit formula for $\langle K \rangle$ via smoothing in an oriented way at each crossing. The result of this smoothing is the set of Seifert circles for the diagram weighted by $s$ or $\overline{s}$ at each smoothing site. Thus the evaluation is given by the formula

$$\langle K \rangle = s^{w(K)} \delta^{SC - 1}.$$

\noindent Substituting, 

$$f_k = (s \delta)^{-w(K)} \langle K \rangle = s^{-w(K)} s^{w(K)} \delta^{-w(K)} \delta^{SC - 1}= \delta^{SC - w(K) - 1}.$$

\noindent Here, $SC$ stands for the number of Seifert circles produced from an oriented link diagram. If $\delta = 1$ then the invariant is trivial. We now show that the invariant is still trivial for $-1$ regardless of the number of Seifert circles or the writhe. To do this we show that the number of Seifert circles minus the writhe of the knot is always odd. Notice that the unknot has one Seifert circle and no writhe, so $f_{unknot} = \delta^{1-0-1} = \delta^{0} = 1$. Similarly, note that the trefoil knot with $w(K) = +3$ produces two Seifert circles, so $f_{trefoil} = \delta^{2-3-1} = \delta^{-2} = +1$. The only question that concerns us is to prove the congruence $$SC - w(K) - 1 \equiv 0\pmod 2,$$ holds for all cases. This is implied by the Lemma in the Appendix to this paper.\\

\noindent Therefore, $SC - w(K) - 1$ is always even and 

$$f_k = \delta^{SC - w(K) - 1} = (\pm 1)^{SC - w(K) - 1} = 1,$$

\noindent thus the invariant is trivial. \qed

\end{proof}

\subsection{\bf The Swap Case}

%\begin{figure}
%	\[
%	\xy
%	(-7,30)*{}; (0,20)*{} **\dir{-}; 
%	(0,30)*{}; (-7,20)*{} **\dir{-}; 
%	(7,30)*{}; (7,20)*{} **\dir{-};
%	(-7,20)*{}; (-7,10)*{} **\dir{-};
%	(0,20)*{}; (7,10)*{} **\dir{-};
%	(7,20)*{}; (0,10)*{} **\dir{-};
%	(-7,10)*{}; (0,0)*{} **\dir{-};
%	(0,10)*{}; (-7,0)*{} **\dir{-};
%	(7,10)*{}; (7,0)*{} **\dir{-};
%	(-7,20)*+{\bullet}="4"; 
%	(0,20)*+{\bullet}="6";
%	(0,10)*+{\bullet}="4"; 
%	(7,10)*+{\bullet}="6";
%	(-7,0)*+{\bullet}="4"; 
%	(0,0)*+{\bullet}="6";
%	(-3,0)*{G};
%	(-10,0)*{F};
%	(3,10)*{F};
%	(10,10)*{G};
%	(-3,21)*{G};
%	(-10,21)*{F};
%	(15,14)*{=};
%	(21,30)*{}; (21,20)*{} **\dir{-};
%	(28,30)*{}; (35,20)*{} **\dir{-}; 
%	(35,30)*{}; (28,20)*{} **\dir{-}; 
%	(21,20)*{}; (28,10)*{} **\dir{-}; 
%	(28,20)*{}; (21,10)*{} **\dir{-}; 
%	(35,20)*{}; (35,10)*{} **\dir{-}; 
%	(21,10)*{}; (21,0)*{} **\dir{-}; 
%	(28,10)*{}; (35,0)*{} **\dir{-}; 
%	(35,10)*{}; (28,0)*{} **\dir{-}; 
%	(28,20)*+{\bullet}="4"; 
%	(35,20)*+{\bullet}="6";
%	(21,10)*+{\bullet}="4"; 
%	(28,10)*+{\bullet}="6";
%	(35,0)*+{\bullet}="4"; 
%	(28,0)*+{\bullet}="6";
%	(31,0)*{F};
%	(38,0)*{G};
%	(24,10)*{G};
%	(17,10)*{F};
%	(31,21)*{F};
%	(38,21)*{G};
%	\endxy
%	\]
%	\caption{The swap relation decomposition of the Yang-Baxter equation. Given that the middle strands must be equivalent to each other, we have the relation $FG=GF$.}
%	\label{fig:swapdecomp}
%\end{figure}

\begin{figure}
     \begin{center}
     \begin{tabular}{c}
     \includegraphics[width=6cm]{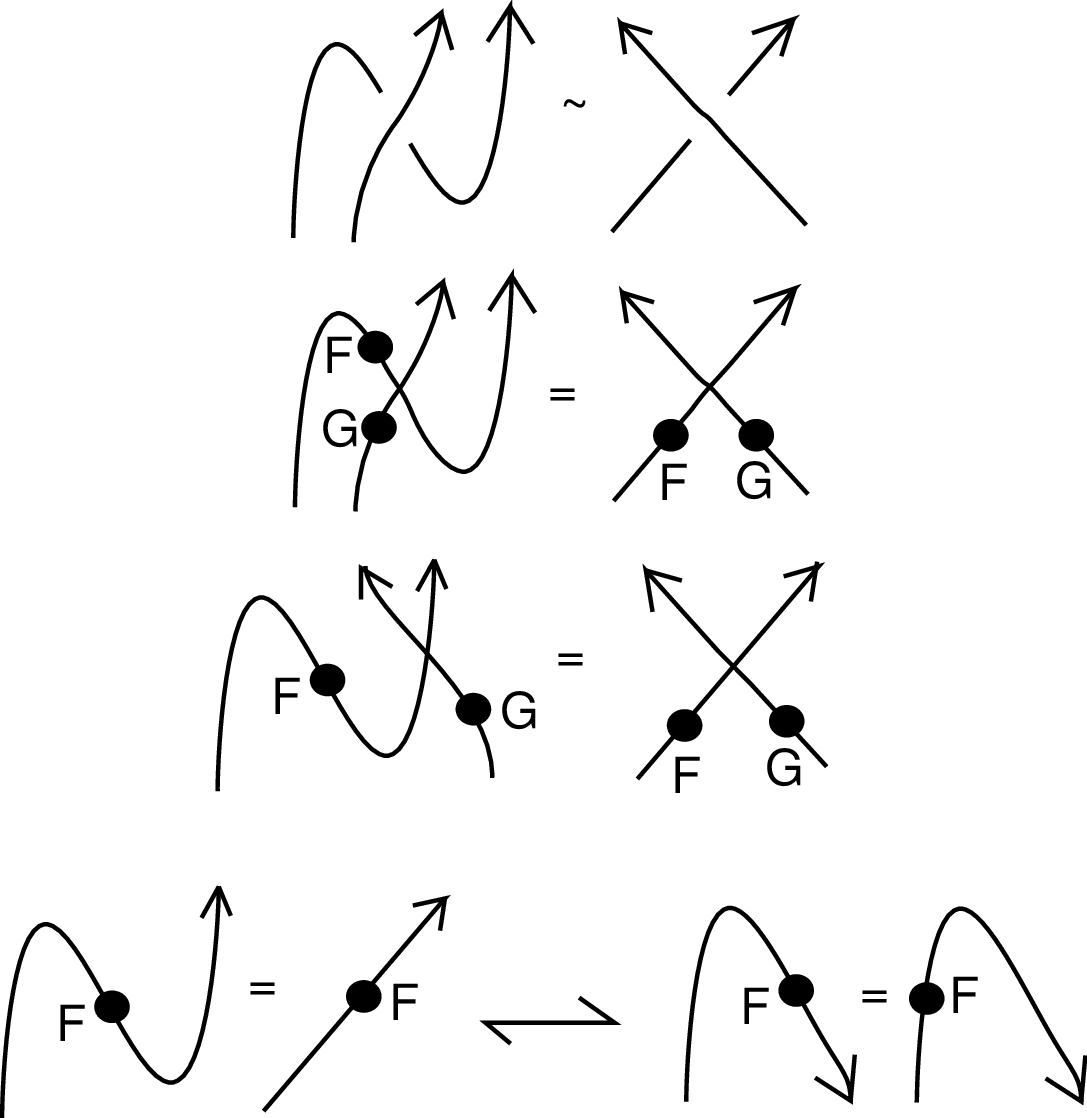}
     \end{tabular}
     \caption{\bf Oriented Swing Move Implies Algebraic Moveability.}
     \label{oswing}
\end{center}
\end{figure}

\textbf{Theorem 4.2} \textit{The state summation model given by $R = (F \otimes G) \circ S$ is trivial for oriented knots.}

\begin{proof}

We use the relationship given in Figure~\ref{fig:swaprelate} and Figure~\ref{yangbaxterswap} except with orientation going upward such that the arrows are on the top endpoints. By using this decomposition of the Yang-Baxter equation, we arrive at the fact that $FG=GF$ (see Figure~\ref{yangbaxterswap}).  View Figure~\ref{oswing} and note that it follows that algebra elements can be moved along the lines of the diagram and collected on a single arc, as we saw in the previous unoriented swap case. (We have indicated one case of a number of cases that are verified in the same way.) We also know that $FG = GF$ and so we can intrepret this freedom by saying that for each crossing we can remove either $FG$ or $\overline{F}\overline{G}$ and collect it directly to an algebra product that we accumulate from the crossing of the diagram. We can remove them from the diagram itself and write that $\langle \raisebox{-0.2\height}{\includegraphics[width=0.5cm, height=0.5cm]{crossorient.eps}} \rangle = FG \langle \raisebox{-0.2\height}{\includegraphics[width=0.5cm, height=0.5cm]{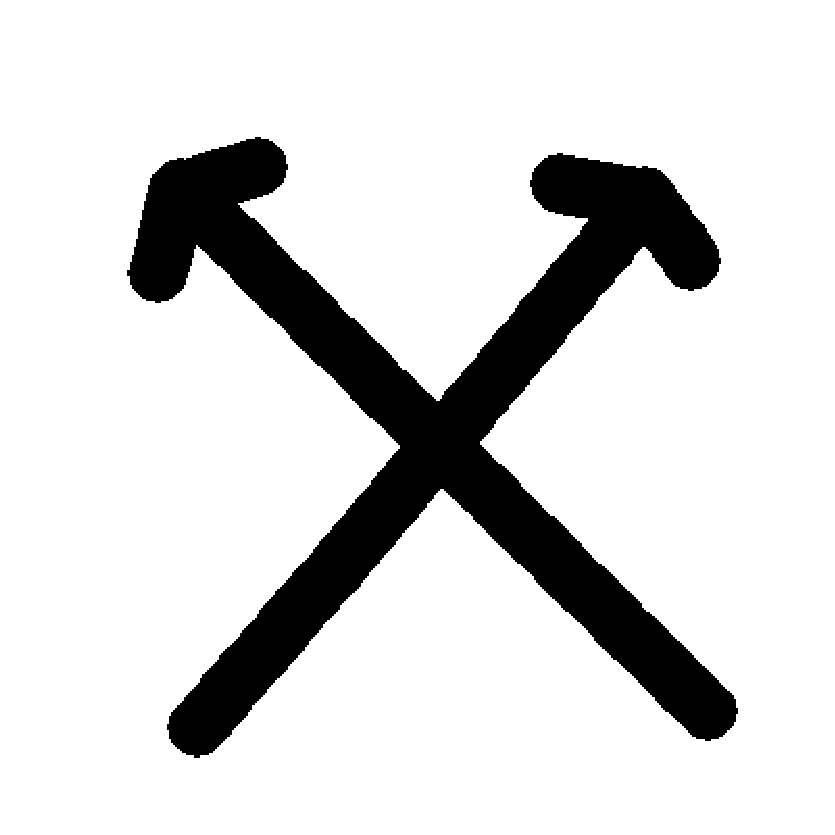}} \rangle$. For each $R$ we then have an $FG$ and for each $\overline{R}$ we have an $\overline{F} \overline{G}$. Note that an oriented $\overline{R}$ has a writhe of $+1$.  Therefore, we have that 

$$\langle \raisebox{-0.2\height}{\includegraphics[width=0.5cm, height=0.5cm]{orientloop.eps}} \rangle = 
\overline{FG} \langle \raisebox{-0.2\height}{\includegraphics[width=0.5cm, height=0.5cm]{lineorient.eps}} \rangle.$$\\

For each positive writhe loop we output a $\overline{FG}$. Our invariant has the form $ f_K = (FG)^{w(K)} (FG)^{N-P}$ where $P$ is the number of oriented $\overline{R}$ crossings and $N$ is the number of oriented $R$ crossings. As all oriented $\overline{R}$ and $R$ crossings have writhe $+1$ and writhe $-1$ respectively, this invariant becomes trivial as both $(FG)^{w(K)}$ and $(FG)^{N-P}$ cancel each other out. \qed

\end{proof}

\section{The Jones Polynomial and Quantum Computation}

We now review \cite{quantumcompute} which gives a local unitary representation that can be used to compute the Jones polynomial \cite{Jones,Jones1,Jones2,Jones3,Jones4} for closures of 3-braids. The quantum computation devolves into finding the trace of a unitary transformation. The result of this construction is a quantum computational model for the Jones polynomial evaluation on a significant class of knots and links that 
is not involved with quantum entanglement since the unitary transformations are in $SU(2).$ This result is very interesting to us, even though it is a special case. We do not know how to obtain the Jones polynomial for all  
knots and links in this way, avoiding entangling operators. The model given here, can be extended to the well-known Fibonacci model \cite{qdeformed}
for quantum computing, but then the transformations are in other unitary groups, and it remains to 
analyse the full role of quantum entanglement in these generalizations.\\ 

The idea behind this construction depends upon the algebra generated by two single qubit density matrices (ket-bras). Let $\ket{v}$ and $\ket{w}$ be two qubits in $V$, a complex vector space of dimension two over the complex numbers. Let $P = \ket{v} \bra{v}$ and $Q = \ket{w} \bra{w}$ be the corresponding ket-bras. Note that as

$$P^{2} = |v|^{2} P,$$
$$Q^{2} = |w|^{2} Q,$$
$$PQP = | \braket{v|w} |^{2} P,$$
$$QPQ = | \braket{v|w} |^{2} Q.$$

\noindent $P$ and $Q$ generate a representation of the Temperley-Lieb algebra. One can adjust parameters to to make a representation of the 3-strand braid group in the form

$$s_1 \rightarrow rP + sI,$$
$$s_2 \rightarrow tQ + uI,$$

\noindent where $I$ is the identity mapping on $V$ and $r,s,t,u$ are suitably chosen scalars. In the following, we use this method to adjust such a representation so that it is unitary. Note also that this is a local unitary representation of $B_3$ to $U(2)$. We leave it as an exercise for the reader to verify that it fits into the general classification of such representations as given in \cite{May}. \\

The representation depends on two symmetric but non-unitary matrices $U_{1}$ and $U_{2}$ with

$$ U_1 = \left[ \begin{array}{cc}
d & 0 \\ 
0 & 0 
\end{array} \right]
= d \ket{w} \bra{w},
$$

\vspace{0.25cm}

$$ U_2 = \left[ \begin{array}{cc}
d^{-1} & \sqrt{1-d^{-2}} \\ 
\sqrt{1-d^{-2}} & d-d^{-1} 
\end{array} \right]
= 
d \ket{v} \bra{v},
$$

\noindent where $w=(1,0)$, and $v = (d^{-1}, \sqrt{1-d^{-2}})$, assuming the entries of $v$ are real. Note that $U^{2}_{1} = d U_1$ and $U^{2}_{2} = d U_1$. Moreover, $U_1 U_2 U_1 = U_1$ and $U_2 U_1 U_2 = U_1$. This is an example of a specific representation of the Temperley-Lieb algebra. The desired representation of the Artin braid group is given on the two braid generators for the 3-strand braid group by the equations:

$$\Phi(s_1) = AI + A^{-1} U_1,$$
$$\Phi(s_2) = AI + A^{-1} U_2,$$

\noindent where $I$ denotes the $2 \times 2$ identity matrix.\\

For any $A$ with $d = -A^2 - A^{-2}$ these formulas define a representation of the braid group. With $A = e^{i \theta}$, we have $d = -2 \cos(2 \theta)$. One finds a specific range of angles $\theta$ in the following disjoint union of angular intervals

$$ \theta \in [0, \pi/6] \sqcup [\pi/3, 2\pi/3] \sqcup [5\pi/6, 7\pi/6] \sqcup [4\pi/3, 5\pi/3] \sqcup [11\pi/6, 2\pi]$$

\noindent that give unitary representations of the 3-strand braid group. Thus, a specialization of a more general representation of the braid group gives rise to a continuous family of unitary representations of the braid group.\\

Note that the traces of these matrices are given by the formulas $Tr(U_1) = Tr(U_2) = d$ while $Tr(U_1 U_2) = Tr(U_2 U_1) = 1$. If $b$ is any braid, let $I(b)$ denote the sum of the exponents in the braid word that expresses $b$. For $b$ a 3-strand braid, it follows that 

$$\Phi(b) = A^{I(b)} I + \Pi(b),$$

\noindent where $I$ is the $2 \times 2$ identity matrix and $\Pi(b)$ is a sum of products in the Temperley-Lieb algebra involving $U_1$ and $U_2$. Since the Temperley-Lieb algebra in this dimension is generated by $I, U_1, U_2, U_1 U_2,$ and $U_2 U_1$, it follows that the value of the bracket polynomial of the closure of the braid $b$, denoted $\langle \overline{b} \rangle$, can be calculated directly from the trace of this representation, except for the part involving the identity matrix. The bracket polynomial evaluation depends upon the loop counts in the states of the closure of the braid, and these loop counts correspond to the traces of the non-identity Temperley-Lieb algebra elements. Note that the closure of the 3-strand diagram for the identity braid in $B_3$ has bracket polynomial $d^2$. The result is the equation 

$$\langle \overline{b} \rangle = A^{I(b)} d^2 + Tr(\Pi(b)),$$

\noindent where $\overline{b}$ denotes the standard braid closure of $b$, and the sharp brackets denote the bracket polynomial. Since the trace of the $2 \times 2$ identity matrix is $2$, we see that 

$$ \langle \overline{b} \rangle = Tr(\Phi(b)) + A^{I(b)} (d^2 - 2).$$

It follows from this calculation that the question of computing the bracket polynomial for the closure of the 3-strand braid $b$ is mathematically equivalent to the problem of computing the trace of the unitary matrix $\Phi(b)$. Therefore, we can define topological invariants from quantum situations that lack any sort of entanglement at all as this calculation depends solely on a single qubit.

By using the method we have described in this section, we show that there is indeed a disparity between topological entanglement and entangling quantum gates. Once we leave the Yang-Baxter formalism it is possible to construct strong topological invariants from non-entangling quantum gates. This phenomena needs further exploration, particularly in regard to the Fibonacci model \cite{qdeformed} and \cite{fibonacci}

\section{Bracket Quantum Link Invariants and Quantum Entanglement}

Recall that our cup and cap matrices are given by the following matrix in the unoriented case:

$$M= \left[ \begin{array}{cc}
0 & iA \\
-iA^{-1} & 0 
\end{array} \right].
$$

\noindent Moreover, the bracket relation in \cite{state} can be given in terms of quantum link invariants as follows

$$
R^{ab}_{cd} = A M^{ab} M_{cd} + A^{-1} \delta^{a}_{c} \delta^{b}_{d}.
$$

\noindent By substituting in our cup/cap matrix and the $2 \times 2$ identity matrix we can give an explicit $R$ as 

$$R= \left[ \begin{array}{cccc}
A^{-1} & 0 & 0 & 0 \\
0 & A^{-1}-A^3 & A & 0 \\
0 & A & 0 & 0 \\
0 & 0 & 0 & A^{-1}
\end{array} \right].
$$

In order for $R$ to be unitary, note that $ (A^{-1} -A^{3})A^{-1} = 0$ or alternatively $1=A^{4}$. Therefore, choosing $A = \pm i$ gives us a unitary, invertible matrix. However, given these choices of $A$ the matrix becomes unentangling as a matrix of the form

$$R = \left[ \begin{array}{cccc}
a\, & 0\, & 0\, & 0\, \\
0 & 0 & d & 0 \\
0 & c & 0 & 0 \\
0 & 0 & 0 & b
\end{array} \right]
$$

\noindent is only entangling when $ab \neq cd$ as shown in \cite{notsoshort}. There is, therefore, no $R$ matrix solution to the bracket that can be an entangling operator. 

\section{Virtual Knot Theory and Quantum Entanglement}

Take the matrix given below:

$$R= \left[ \begin{array}{cccc}
0 & 0 & 0 & A \\
0 & A^{-1} & 0 & 0 \\
0 & 0 & A^{-1} & 0 \\
A & 0 & 0 & 0
\end{array} \right],
$$

\noindent where $A \in S^1$ and $R$ is unitary and a solution to the Yang-Baxter equation. We now show that $R$ is an entangling matrix. Take a decomposed state$$\ket{\psi} = (x\ket{0} + y\ket{1}) \otimes (z\ket{0} + w\ket{1}) = xz\ket{00} + yz\ket{10} + xw\ket{01} + yw\ket{11}.$$
\noindent Now, we apply $R$ to $\ket{\psi}$ to get

$$
R\ket{\psi} = xzA\ket{11} + yz A^{-1}\ket{10} + xw A^{-1}\ket{01} + yw A \ket{00}.
$$

\noindent From our definition of entanglement, we take the determinant of the resultant state
$$\det \left[ \begin{array}{cc}
A yw & A^{-1} xw \\
yz A^{-1} & xz A 
\end{array} \right]
= 
xyzw(A^2 - A^{-2}).
$$

\noindent We must have that $A = e^{i \theta}$, which implies that 
\begin{center}
$$
\begin{array}{rcl} A^2 - A^{-2} & = & e^{2i \theta} - e^{-2i \theta} \\ & = & \cos(2 \theta) + i \sin(2\theta) - \cos(2 \theta) + i \sin(2 \theta) \\ & = & 2i \sin(2 \theta), \end{array}
$$

\end{center}

\noindent which shows that $\sin(2 \theta) \neq 0$. Therefore, there are a continuum of solutions such that this given $R$ matrix is entangling. In \cite{virtual}, $R$ is shown to only detect the writhe of classical knots; however, when $R$ is applied to virtual knots it is a much stronger invariant. Many of the arguments of this paper can be generalized to virtual knots. Moreover, the relationship between physics and virtual knots has yet to be explored in detail. For quantum computing, the virtual crossing can be modeled as a swap gate (interchange qubits as in $S \ket{01} = \ket{10})$. Thus it is natural to use the virtual braid group and its unitary representations for quantum computing. We will return to this subject in a subsequent paper.

\section{Topological Entanglement and Quantum Entanglement}
In this paper we have, so far, discussed topological entanglement and 
quantum entanglement by examining quantum operators that are solutions to the Yang-Baxter equation. The operators $R$ that we have considered are unitary solutions to the 
Yang-Baxter equation that act of the tensor product of a single qubit space $V$ with itself. Such an operator $R: V\otimes V \longrightarrow  V\otimes V$ can be an entangling operator 
in the quantum sense. We have shown in this paper that such operators will not produce non-trivial invariants of knots and links (except in very special cases) unless they are quantum entangling. This establishes a connection between the ability to detect topological entanglement and entangled quantum states.\\

In this section of the paper we discuss more generally the theme of quantum entanglement and topological entanglement. 
We begin with the Aravind hypothesis \cite{Ara}.
The Aravind hypothesis suggests that topological linking may be directly comparable to quantum entanglement. We discuss the pros and cons of this hypothesis below.
The main work of this paper shows that there is a relationship between quantum entangling operators and invariants of knots and links. The Aravind hypothesis suggests that 
there may be a more direct relationhip of topological and quantum entanglement.\\

We then discuss the relationship of space, spacetime and quantum entanglement in the context of the a hypothesis of Susskind and Malcedena \cite{ER}. This  $ER = EPR$ hypothesis is based on the suggestion that the connectivity of spacetime
is a phenomenon of quantum entanglement. Susskind asserts that the entanglement of distant particles is equivalent to the existence of an Einstein-Rosen bridge connecting them. If this hypothesis is true, then there is indeed a topological underpinning for quantum entanglement. Here we will only make foundational comments on the $ER = EPR$ hypothesis. In the discussion below we examine entanglement and teleportation in relation to the construction of a that is space augmented by quantum states. Since an entangled state
such as $|\delta \rangle = \frac{1}{\sqrt{2}}(|01\rangle + |10\rangle)$ is formulated without any background space, we point out that it is possible graphically to form a new space from the given space or spactime $S$ of the physics by attaching a corresponding quantum network to $S.$ The new space $S'$ has connectivity related to the entanglement. This construction can then be considered as a precursor to the spacetime with an Einstein-Rosen bridge connecting the sites of the entangled particles.\\

\subsection{\bf The Aravind Hyppothesis}

Link diagrams can be used as holders of information. Aravind \cite{Ara} proposed that
the topological entanglement of a link should correspond to 
the quantum entanglement of a state. Each link component would correspond to a tensor factor of the state. {\em Measurement of a link 
would be  modeled by deleting one
component of the link.} 
A key example is the 
Borommean rings. See Figure~\ref{boro}.
Deleting any component of the Boromean rings 
yields a remaining pair of unlinked rings. The three Borromean 
rings are entangled, but any two of them are unentangled.
In this sense the Borromean rings are analogous to 
the $GHZ$ state $|GHZ\rangle  = (1/\sqrt{2})(|000\rangle  + |111\rangle )$.
Measurement in any factor of the $GHZ$ yields 
an unentangled state. Aravind points out that 
this property is basis dependent. Kauffman and Lomonaco pointed out \cite{notsoshort} 
that there are states whose entanglement 
after measurement is a matter of probability (via quantum amplitudes).
Consider for example the state 

{\bf \[
|\psi \rangle=|001\rangle  + |010\rangle  + |100\rangle .
\]} 

\noindent Measurement in any coordinate yields probabilistically an entangled or an
unentangled state. For example

{\bf \[
|\psi \rangle= |0\rangle(|01\rangle  + |10\rangle)  + |1\rangle |00\rangle .
\] }

\noindent so that projecting to $|1\rangle$ in the first coordinate yields an unentangled state, while projecting to $|0\rangle$ yields an entangled state.
\bigbreak
 
New ways to use link diagrams must be invented to map the properties 
of such states.  One direction is to consider appropriate notions of quantum knots so that one can formulate
superpositions of topological types as in \cite{QK1}. But one needs to go deeper in this consideration.\\

The relationship of topology and physics needs to be examined carefully.
Topological properties of systems are properties
that remain invariant under certain transformations that are identified as ``topological
equivalences". In making quantum physical models, these equivalences should correspond to unitary 
transformations of an appropriate Hilbert space. Accordingly, Kauffman and Lomonaco formulated a model for 
{\it quantum knots} \cite{QK2,QK3,QK4,QK5} that meets these requirements.
A quantum knot system represents the \textquotedblleft quantum embodiment"
of a closed knotted physical piece of rope. \ A quantum knot (i.e., an
element $\left\vert K\right\rangle $ lying in an appropriate Hilbert space $H_{n}$, 
as a state of this system, represents the state of such a knotted
closed piece of rope, i.e., the particular spatial configuration of the knot
tied in the rope. Associated with a quantum knot system is a group of
unitary transformations $A_{n}$, called the {\it ambient group},
which represents all possible ways of moving the rope around (without
cutting the rope, and without letting the rope pass through itself.) Unlike a classical closed piece of rope, a quantum knot can exhibit
non-classical behavior, such as quantum superposition and quantum
entanglement.  
The {\it knot type} of a quantum knot $\left\vert K\right\rangle $ is the orbit of the quantum knot under the action of the ambient group 
$A_{n}$. This leads to new questions connecting quantum
computing and knot theory. \\

\begin{figure}
     \begin{center}
     \begin{tabular}{c}
     \includegraphics[height=4cm]{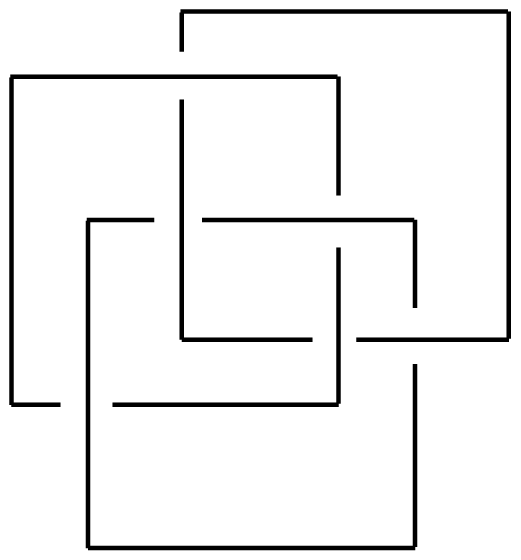}
     \end{tabular}
     \end{center}
     \caption{\bf Borromean Rings }
     \label{boro}
     \end{figure} 
     \bigbreak

\begin{figure}
     \begin{center}
     \begin{tabular}{c}
     \includegraphics[height=4cm]{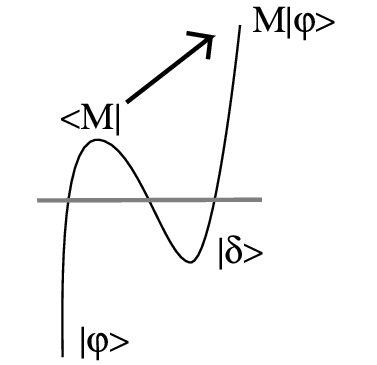}
     \end{tabular}
     \end{center}
     \caption{\bf A Teleportation Scenario}
     \label{teleport}
     \end{figure} 
     \bigbreak

\begin{figure}
     \begin{center}
     \begin{tabular}{c}
     \includegraphics[height=4cm]{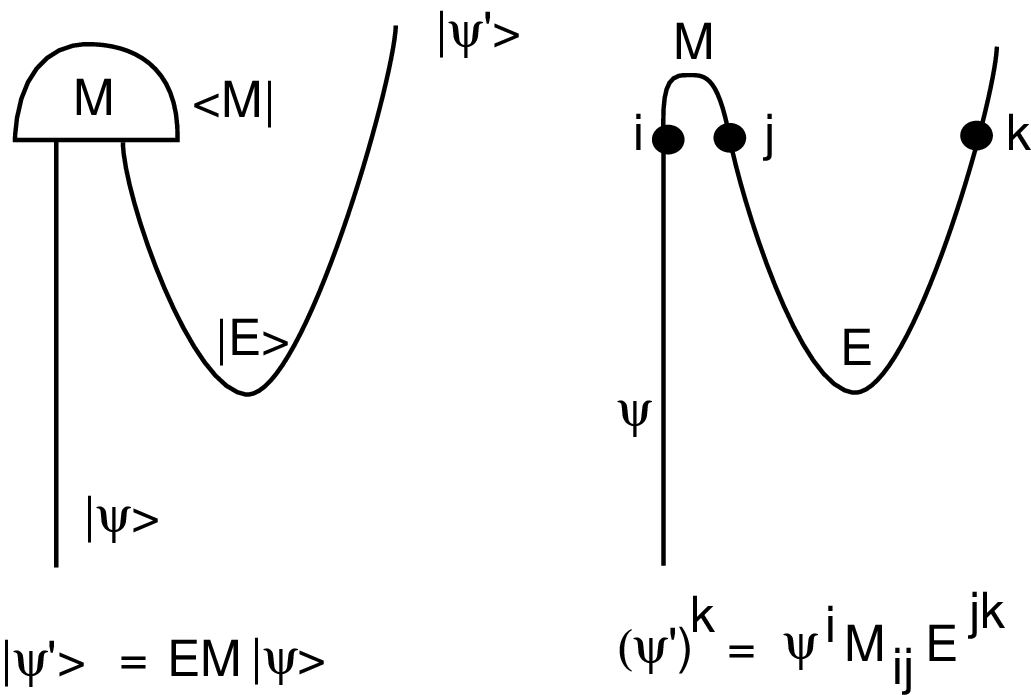}
     \end{tabular}
     \end{center}
     \caption{\bf Teleportation Tensors}
     \label{teleportens}
     \end{figure} 
     \bigbreak

\begin{figure}
     \begin{center}
     \begin{tabular}{c}
     \includegraphics[height=4cm]{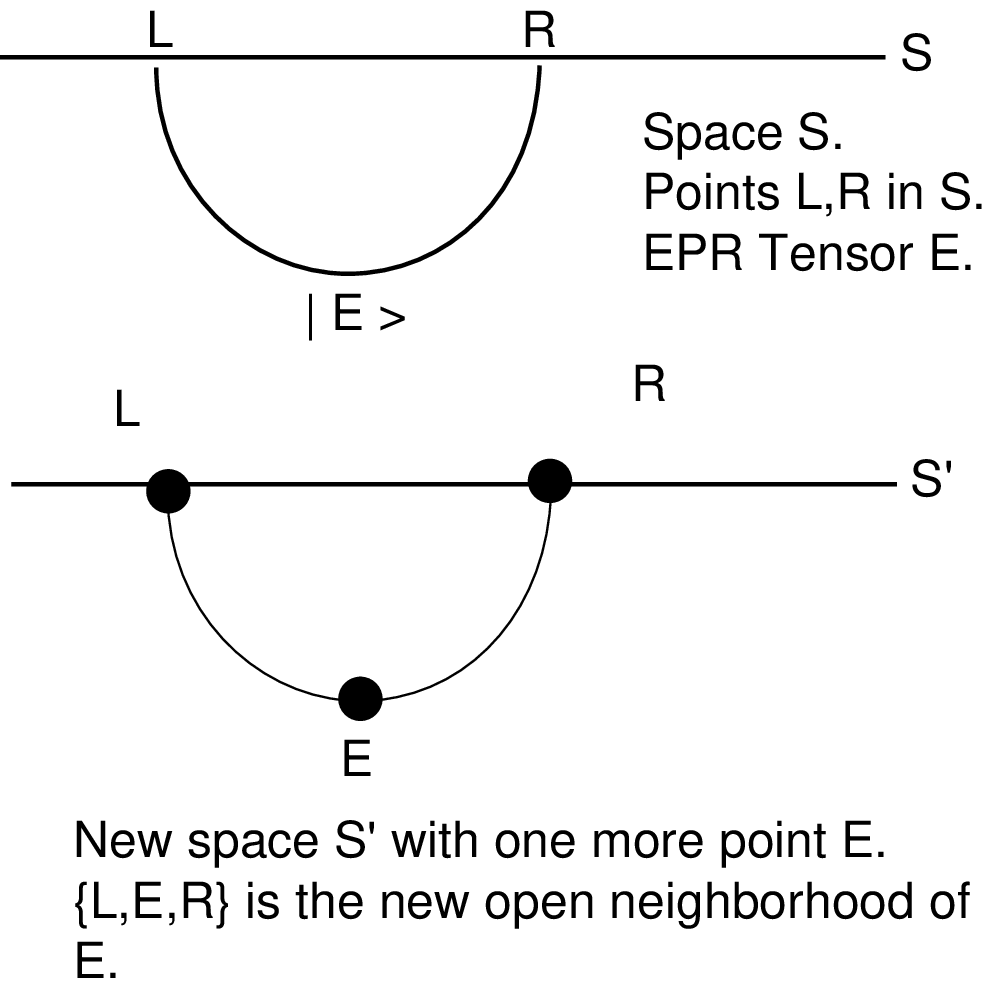}
     \end{tabular}
     \end{center}
     \caption{\bf Augmented Space}
     \label{augment}
     \end{figure} 
     \bigbreak

\begin{figure}
     \begin{center}
     \begin{tabular}{c}
     \includegraphics[height=4cm]{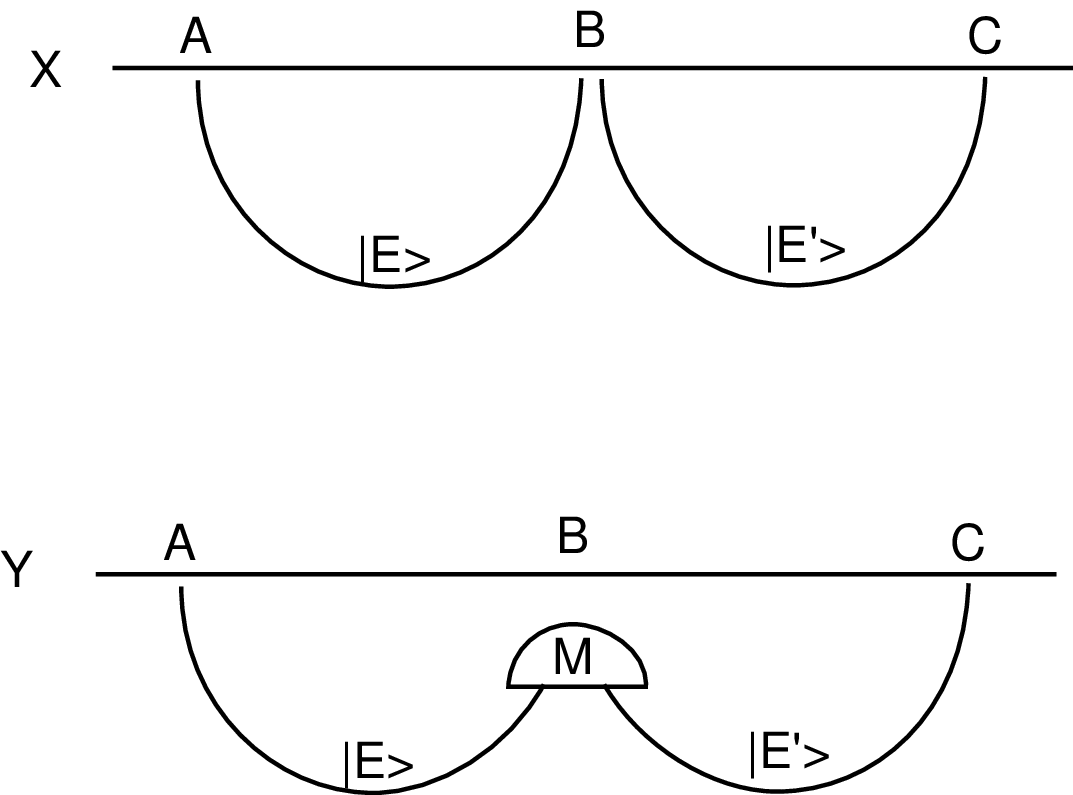}
     \end{tabular}
     \end{center}
     \caption{\bf Entanglement Swap}
     \label{swap}
     \end{figure} 
     \bigbreak

\subsection{\bf Space, Time, Quantum Networks and Entanglement}
Here is a summary for understanding quantum teleportation. Take  $$|\delta \rangle = \frac{1}{\sqrt{2}}(|01\rangle + |10\rangle)$$ as a representative 
entangled state. Regard $|\delta\rangle$ as representing the state of two particles that we shall call $L$ (left) and $R$ (right) corresponding to $\delta$'s right and left tensor factors.
Measuring $|\delta\rangle$ results either in $|01\rangle$ or $|10\rangle.$ If an observer measures the left particle, and sees $0$ then an observer who will measure the right particle
must see $1$ and vice versa. Nowhere in the quantum state $|\delta\rangle$ is there any information about the distance between the particles $L$ and $R$ or any information about the relative times for measurements to occur at the locales for these particles.\\

Note that the entangled state $|\delta\rangle$ is in the tensor product $V \otimes V$ where $V$ is a qubit space spanned by $|0\rangle$ and $|1\rangle.$ A general element in
$V \otimes V$ has the form $$| A \rangle = a_{00}|00\rangle + a_{01}|01\rangle + a_{10}|10\rangle + a_{11}|11\rangle,$$ and can presented as a $2 \times 2$ matrix
$$ A = \left[ \begin{array}{cc}
a_{00} & a_{01} \\
a_{10} & a_{11} 
\end{array} \right].$$
Thus the matrix for $|\delta\rangle$ is the identity matrix
$$ I = \left[ \begin{array}{cc}
1 & 0\\
0 &  1 
\end{array} \right].$$
By the same token, a successful measurement on two tensor lines can be represented in the dual basis spanned by elementary bras as 
$$\langle M | = m_{00} \langle 00 | + m_{01} \langle 01 |  + m_{10} \langle 10 | + m_{11} \langle 11|,$$ with corresponding matrix
$$ M= \left[ \begin{array}{cc}
m_{00} & m_{01} \\
m_{10} & m_{11} 
\end{array} \right].$$
Now consider Figure~\ref{teleport} where we have indicated an initial qubit state $|\phi \rangle $ tensored with the entangled state $| \delta \rangle.$
A sucessful measurement has been made on the first two tensor lines. We assert that the state on the final tensor line is given by $M |\phi \rangle$ where this denotes the 
action of the matrix $M$ of the measurement $\langle M |$ on the vector $| \phi \rangle.$ The reader will find the details of this calcluation in \cite{tele}. This means that if Alice is at the site of the left particle and performs the measurement $\langle M |,$ then she knows that Bob (at the site of the right particle) will have the quantum state  $M |\phi \rangle.$ If the matrix $M$ is invertible and unitary, Alice can phone Bob and tell him to apply $M^{-1}$ to the state that he has. The result will be that Bob will then have a perfect copy of the original state $| \phi \rangle.$ This is the key to teleportation. Specific teleportation protocols use an orthonormal  measurement basis (for example the so-called Bell basis) such that all the matrices of the basis elements are unitary. Then this teleportation protocal can be applied whenever Alice measures, at her end, the two left tensor lines. It is as if the wiggle in the line in Figure~\ref{teleport} is pulled staight and the new straightened line represents the transformation $M.$ There is a geometry in the tensor 
diagrams for the teleportation procedure. It is this geometry that we wish to pursue to understand the geometry and topology of entanglement.\\

In Figure~\ref{teleportens} we illustrate the general case for a single qubit teleportation. The entangled state now has matrix $E,$ not neccessarily the identity, an the measurement has matrix $M.$ We then see from the figure that $|\psi '\rangle = EM|\psi \rangle.$  The matrix  $E$ of an entangled state is neccessarily invertible, and so when $E$ and $M$ are
unitary, our previous description of the teleportation procedure goes over mutatis mutandis. 
Using indices, the description of the state transformation is given by the equation 
$$(\psi ')^{k} = \psi^{i} M_{ij}E^{jk}.$$ The important point to note about this index version of the equation is that it is an exact translation of the structure of the tensor network
given on the right part of the figure. This tensor network is the detailed expression of the tensor diagram on the left part of Figure~\ref{teleportens}. This transformation from 
 $(\psi^{i})$ to $(\psi ')^{k} = \psi^{i} M_{ij}E^{jk}$ can be described by following the connectivity of the tensor network from Alice's locale to Bob's locale. The successful measurement
 $\langle M |$ completes the connection and transforms the quantum information at $|\psi \rangle$, located with Alice to $|\psi '\rangle = EM|\psi \rangle,$ located with Bob.
 It is a transfer of quantum information, a transfer of quantum states. To obtain observed information transfer one would need to control both measurement at Alice's end and corresponding measurement at Bob's end. Nevertheless, the tensor network for the entanglement can be viewed as a way to augment the simple space between Alice and Bob. This extra connectivity between Alice and Bob resides in the entangled state $|E\rangle$ that connects them.\\
 
 We formalize the idea that the tensor network for quantum entanglement can augment the original physical space to create a new connectivity. Let $S$ be the given background space for the physical locations of particles. For each entangled state $|E\rangle$ with corresponding particles (approximately) located at points $L$ and $R$ in the space $S$  associate a new point $E$ and a new open neighborhood $\{ L, E, R \}$ for this new point $E.$ Let $S'$ be the new space with topology generated by these new neighborhoods of points corresponding to the entangled states. We call $S'$ the {\it quantum tensor space} associated with $S$ and its quantum network.
See Figure~\ref{augment} for an illustration of this concept.\\

To construct the quantum tensor space, we introduce a least topological structure that can produce the special connection between $L$ annd $R.$ Note that the new space has non-Hausdorff points for each entangled state. The neighborhood  $\{ L, E, R \}$ is a combnatorial topological analogue of an Einstein-Rosen bridge connecting $L$ and $R.$ The analogy is important. Note that an observer in the space $S',$ cannot move continously from $L$ to $E$ without invoking an open neighborhood of $E$ and the least such neighborhood contains $R.$ Letting Alice be the observer at $L$ and Bob the observer at $R,$ we can say that Alice and Bob can meet together at the connecting point $E$ in the analogue black hole. The point $E$ is the analogue of the event horizon of an Einstein-Rosen bridge between $L$ and $R.$ We will explore the analogies between connectivity in the the quantum tensor spaces and connectivity via Einstein-Rosen bridges in a later paper. It is possible that for larger networks and states with many particles these precursors to Einstein-Rosen Bridges will approximate the bridges in the continuum spacetime. For our purposes, we introduce this formalism to show how it is possible to weld a combinatorial quantum tensor network to a given background space.\\

Figure~\ref{swap} illustrates the procedure known as {\it entanglement swapping}. Locations $A$ and $B$ are connected by a entangled state $|E\rangle$ and 
locations $B$ and $C$ are connected by an entangled state $|E'\rangle.$ By performing a measurement $\langle M|$ at $B$ we connect the two entangled states and make
a new entangled state that connects $A$ with $C.$ In the process, the entanglement connection with $B$ is lost. This example shows how the topology of the quantum tensor space will change under the act of measurement. Just as the quantum network undergoes graphical cut and rejoin operations under measurement, the corresponding quantum space, made by the prescription above, will change its connectivity properties. The actions on our simple spaces are easy to understand. In the case of the $ER = EPR$ hypothesis it will be very interesting to 
see what is the meaning of a procedure such as entanglement swapping.\\
 
 The purpose of this section of the paper has been to show how topological connectivity in the form of knots and in the form of spatially realized quantum networks can be used to make new points of view about quantum entanglement. These points of view shed new light on both the Aravind hypothesis and the Susskind $ER = EPR$ hypothesis. In our opinion the Susskind hypothesis is the deeper of the two and most likely to lead to new physics. Nevertheless, the connections between knot theory and the structure of quantum information
 are very strong and deserve further investigation.\\

\section{Summary}

We have shown that entanglement is a necessary condition for forming invariants from $R$ matrices from state summation models in the oriented case, while the arguments used by \cite{amsshort} do not generalize to the unoriented case for state summation models. We must highlight the fact that this appears to be the first time that the two methods (combinatorial and enhanced Yang-Baxter operators) have been shown to differ. We also have found that there is a potential relationship between virtual knots and quantum entanglement that could elucidate more about the relationship between topology and quantum entanglement. However, there exist quantum algorithms for forming topological invariants of knots that rely on no entanglement at all, as in the quantum algorithm described here for computing the Jones polynomial on three-strand braids, which depends only on a single qubit. In conclusion, by studying the boundary between topological and quantum entanglement we can construct a correspondence between topological invariants and entangling $R$ matrices that may have a significant impact on the study of quantum computing. In the final section of the paper we have discussed relations between the ideas of this paper and the entanglement hypotheses of Aravind and the ER = EPR hypothesis of Susskind and his collaborators. In the light of the latter hypothesis we have shown how to augment a space to a new space that contains a topological version of the tensor networks describing its quantum structure.\\

\section{Acknowledgments}  Eshan Mehrotra takes pleasure in thanking his teachers, mentors, and parents. We would also like to thank Dr. Vandana Chinwalla and Dr. Robyn Fischer for their part in facilitating the Student Inquiry and Research program at the Illinois Math and Science Academy. This work could not have been possible without the environment of free inquiry they labored to foster. It gives Louis Kauffman great pleasure to thank Stathis Antoniou for stimulating conversations in the course of this research.
Kauffman's work was supported by the Laboratory of Topology and Dynamics, Novosibirsk State University (contract no. 14.Y26.31.0025 with the Ministry of Education and Science of the Russian Federation).

\section{APPENDIX}
Recall that $w(K)$ is the writhe of an oriented diagram $K$ and that $SC(K)$ denotes the number of Seifert circuits obtained from the diagram $K$ by smoothing all crossings in an oriented manner, and $cr(K)$ is the number of crossings in the diagram $K.$ Furthermore, $P$ and $N$ denote the number of positive and negative crossings of $K.$

\noindent We shall prove the following\\ 

\noindent {\bf Lemma.}  For an oriented link diagram $K,$ $$SC(K) - w(K) - 1 \equiv 0\pmod 2.$$

\noindent {\it Proof} The above equation implies that $$SC(K) \equiv w(K) + 1\pmod 2.$$ Note that the parity of the writhe of a knot and its crossing number are the same, as $$w(K) = P - N \equiv P + N = cr(K)\pmod 2.$$ Thus, we wish to prove that $SC(K)  \equiv cr(K) + 1.$ From now on, we write $SC$ for $SC(K).$\\

Consider a knot diagram where the crossings have been replaced with flat nodes. Note that by Euler's theorem we have that $$v+2=R,$$ where $R$ is the number of regions in the diagram and $v$ is the number of nodes ($cr(K) = v$). This theorem implies that $$v \equiv R\pmod 2.$$ 
\noindent It follows  that $$SC \equiv R+1\pmod 2.$$ For the unknot, $R = 2$ and $SC = 1$, as shown below. 

\[
\xy
(0,0)*\xycircle(7,7){-};
(15,2)*{R=2};
(15,-2)*{SC = 1};
\endxy 
\]

\noindent Imagine that each of the following diagrams is a subsection of a much larger diagram. Take the oriented loop,

\[
\xy
{\ar@/^-1pc/(5,-5)*{};(5,5)*{}};
{\ar (5,5)*{}; (-5,-5)*{}}; 
{\ar (-5,5)*{}; (5,-5)*{}}; 
\endxy
\]

\noindent Suppose that in this diagram we had that $SC \equiv R + 1\pmod 2$. After we apply our smoothing we get the following diagram. 

\[
\xy
{\ar (-5,5)*{}; (-5,-5)*{}}; 
{\ar@/^-1pc/(5,-5)*{};(5,5)*{}};
{\ar@/^-1pc/(5,5)*{};(5,-5)*{}};
\endxy
\]

\noindent Note that this diagram has $SC+1$ Seifert circles and $R+2$ regions, so $$SC+1 \equiv (R+1) + 1\pmod 2$$. 

\noindent We now assign the following diagram the value of $R$ and $SC$. We want to show that the fact that $R-SC \equiv 1\pmod 2$ is true does not change under regular isotopy. 

\[
\xy
{\ar@/^1pc/(-5,10)*{};(-5,-5)*{}}; 
{\ar@/^-1pc/(5,10)*{};(5,-5)*{}};
(10,4)*{R};
(10,0)*{SC};
\endxy
\]

\noindent The next diagram represents the oriented Reidemeister II move with the orientation of each strand pointing in the same direction. It has been given values $R'$ and $SC'$. Note that $R' = R +2$ and $SC' = SC$ in this diagram. Therefore, $R' - SC' = R - SC + 2$, so $R'-SC' \equiv 1\pmod 2$. 

\[
\xy
{\ar@/^1pc/(-3,10)*{};(-3,-5)*{}}; 
{\ar@/^-1pc/(3,10)*{};(3,-5)*{}};
(10,4)*{R'};
(10,0)*{SC'};
\endxy
\]

\noindent The next diagram is the oriented Reidemeister II move with the orientations in the opposite direction. 

\[
\xy
{\ar@/^1pc/(-3,10)*{};(-3,-5)*{}}; 
{\ar@/^1pc/(3,-5)*{};(3,10)*{}};
(10,4)*{R'};
(10,0)*{SC'};
\endxy
\]

\noindent For this diagram, we must show that two separate cases hold true. In the first case we must check that the relation remains the same if the bottom left strand connects to the top left and the top right strand connects to the bottom right (so $R' = 3$ and $SC' = 2$). It is easy to see that this case results in $R' = 2$ and $SC' = 1$. In the second case, we connect the bottom left strand to the bottom right and the top right strand to the top left ($R' = 2$ and $SC' = 1$). This case results in $R' = 3$ and $SC' = 2$. In all of these cases we have that $$R' - SC' \equiv 1 \equiv R - S\pmod 2.$$

\noindent We now present the Reidemeister III move. The strands are labeled in the diagram below. There are several cases for this move as shown below.

\[
\xy
{\ar (5,5)*{}; (-5,-5)*{}}; 
{\ar (5,-5)*{}; (-5,5)*{}}; 
{\ar (-7,2)*{}; (7,2)*{}}; 
(14,0)*{\leftrightarrow};
{\ar (29,5)*{}; (19,-5)*{}}; 
{\ar (29,-5)*{}; (19,5)*{}}; 
{\ar (17,-2)*{}; (31,-2)*{}};
(7,7)*{2};
(-7,-7)*{5};
(7,-7)*{4};
(-7,7)*{1};
(-9,2)*{6};
(9,2)*{3};
\endxy
\]

\noindent \textbf{Case 1:} Suppose that 1 connects to 2, 3 connects to 4, and 5 connects to 6. We use the notation $(12)(34)(56)$ to represent the connections. The resulting diagrams for each side of the relation would be 

\[
\xy
{\ar@/^-1pc/(-10,0)*{};(-5,5)*{}};
{\ar@/^-1pc/(5,5)*{};(10,0)*{}};  
{\ar@/^-1pc/(5,-15)*{};(-5,-15)*{}}; 
{\ar@/^-1pc/(-3,-4)*{};(3,-4)*{}}; 
{\ar@/^-1pc/(3,-4)*{};(-3,-4)*{}}; 
{\ar@/^-1pc/(-3,-4)*{};(3,-4)*{}}; 
{\ar@/^1pc/(10,0)*{};(5,-15)*{}};
{\ar@/^1pc/(-5,-15)*{};(-10,0)*{}}; 
{\ar@/^1pc/(-5,5)*{};(5,5)*{}}; 
(18, -4)*{\leftrightarrow};
{\ar@/^1pc/(35,5)*{};(25,5)*{}};
{\ar@/^1pc/(20,-10)*{};(25,-15)*{}};
{\ar@/^1pc/(35,-15)*{};(40,-10)*{}}; 
{\ar@/^1pc/(27,-5)*{};(33,-5)*{}}; 
{\ar@/^1pc/(33,-5)*{};(27,-5)*{}}; 
{\ar@/^1pc/(25,5)*{};(35,5)*{}}; 
{\ar@/^1pc/(25,-15)*{};(20,-10)*{}}; 
{\ar@/^1pc/(40,-10)*{};(35,-15)*{}}; 
\endxy
\]

\noindent In this case, on the left side $R' = 3$ and $SC' = 2$ so $R' - SC' \equiv 1\pmod 2$. For the right side, $R' = 5$ and $SC' = 4$ so $R' - SC' \equiv 1\pmod 2$ too. 

\noindent \textbf{Case 2:} $(12)(36)(45)$ 

\[
\xy
{\ar@/^-1pc/(-10,0)*{};(-5,5)*{}};
{\ar@/^-1pc/(5,5)*{};(10,0)*{}};  
{\ar@/^-1pc/(5,-15)*{};(-5,-15)*{}}; 
{\ar@/^-1pc/(-3,-4)*{};(3,-4)*{}}; 
{\ar@/^-1pc/(3,-4)*{};(-3,-4)*{}}; 
{\ar@/^-1pc/(-3,-4)*{};(3,-4)*{}}; 
{\ar@/^1pc/(-5,5)*{};(5,5)*{}};
{\ar@/^8pc/(10,0)*{};(-10,0)*{}};
{\ar@/^-1pc/(-5,-15)*{};(5,-15)*{}};
(18, -4)*{\leftrightarrow};
{\ar@/^1pc/(35,5)*{};(25,5)*{}};
{\ar@/^1pc/(20,-10)*{};(25,-15)*{}};
{\ar@/^1pc/(35,-15)*{};(40,-10)*{}}; 
{\ar@/^1pc/(27,-5)*{};(33,-5)*{}}; 
{\ar@/^1pc/(33,-5)*{};(27,-5)*{}}; 
{\ar@/^1pc/(25,5)*{};(35,5)*{}};
{\ar@/^3pc/(40,-10)*{};(20,-10)*{}};
{\ar@/^-1pc/(25,-15)*{};(35,-15)*{}};
\endxy
\]

\noindent In this case, on the left $R' = 4$ and $SC' = 3$ and on the right $R' = 4$ and $SC' = 3$. Both sides satisfy the relation. 

\noindent \textbf{Case 3:} $(16)(32)(45)$ 

\[
\xy
{\ar@/^-1pc/(-10,0)*{};(-5,5)*{}};
{\ar@/^-1pc/(5,5)*{};(10,0)*{}};  
{\ar@/^-1pc/(5,-15)*{};(-5,-15)*{}}; 
{\ar@/^-1pc/(-3,-4)*{};(3,-4)*{}}; 
{\ar@/^-1pc/(3,-4)*{};(-3,-4)*{}}; 
{\ar@/^-1pc/(-5,-15)*{};(5,-15)*{}};
{\ar@/^-1pc/(10,0)*{};(5,5)*{}};
{\ar@/^-1pc/(-5,5)*{};(-10,0)*{}};
(15, -4)*{\leftrightarrow};
{\ar@/^1pc/(35,5)*{};(25,5)*{}};
{\ar@/^1pc/(20,-10)*{};(25,-15)*{}};
{\ar@/^1pc/(35,-15)*{};(40,-10)*{}}; 
{\ar@/^1pc/(27,-5)*{};(33,-5)*{}}; 
{\ar@/^1pc/(33,-5)*{};(27,-5)*{}}; 
{\ar@/^-1pc/(25,5)*{};(20,-10)*{}};
{\ar@/^-1pc/(40,-10)*{};(35,5)*{}};
{\ar@/^-1pc/(25,-15)*{};(35,-15)*{}};
\endxy
\]

\noindent The left side has $R' = 5$ and  $SC' = 4$ and the right side has $R' = 3$ and $SC' = 2$. 

\noindent \textbf{Case 4:} $(14)(32)(56)$ 

\[
\xy
{\ar@/^-1pc/(-10,0)*{};(-5,5)*{}};
{\ar@/^-1pc/(5,5)*{};(10,0)*{}};  
{\ar@/^-1pc/(5,-15)*{};(-5,-15)*{}}; 
{\ar@/^-1pc/(-3,-4)*{};(3,-4)*{}}; 
{\ar@/^-1pc/(3,-4)*{};(-3,-4)*{}}; 
{\ar@/^-1pc/(10,0)*{};(5,5)*{}};
{\ar@/^1pc/(-5,-15)*{};(-10,0)*{}};
{\ar@/^-6pc/(-5,5)*{};(5,-15)*{}};
(15, -4)*{\leftrightarrow};
{\ar@/^1pc/(40,5)*{};(30,5)*{}};
{\ar@/^1pc/(25,-10)*{};(30,-15)*{}};
{\ar@/^1pc/(40,-15)*{};(45,-10)*{}}; 
{\ar@/^1pc/(32,-5)*{};(38,-5)*{}}; 
{\ar@/^1pc/(38,-5)*{};(32,-5)*{}}; 
{\ar@/^1pc/(30,-15)*{};(25,-10)*{}};
{\ar@/^-5pc/(30,5)*{};(40,-15)*{}};
{\ar@/^-1pc/(45,-10)*{};(40,5)*{}};
\endxy
\]

\vspace{0.25cm}

\noindent The left side has $R' = 4$ and $SC' = 3$ and the right side has $R' = 4$ and $SC' = 3$

\noindent \textbf{Case 5:} $(52)(16)(34)$ 

\[
\xy
{\ar@/^-1pc/(-10,0)*{};(-5,5)*{}};
{\ar@/^-1pc/(5,5)*{};(10,0)*{}};  
{\ar@/^-1pc/(5,-15)*{};(-5,-15)*{}}; 
{\ar@/^-1pc/(-3,-4)*{};(3,-4)*{}}; 
{\ar@/^-1pc/(3,-4)*{};(-3,-4)*{}}; 
{\ar@/^6pc/(-5,-15)*{};(5,5)*{}};
{\ar@/^1pc/(10,0)*{};(5,-15)*{}};
{\ar@/^-1pc/(-5,5)*{};(-10,0)*{}};
(15, -4)*{\leftrightarrow};
{\ar@/^1pc/(45,5)*{};(35,5)*{}};
{\ar@/^1pc/(30,-10)*{};(35,-15)*{}};
{\ar@/^1pc/(45,-15)*{};(50,-10)*{}}; 
{\ar@/^1pc/(37,-5)*{};(43,-5)*{}}; 
{\ar@/^1pc/(43,-5)*{};(37,-5)*{}}; 
{\ar@/^-1pc/(35,5)*{};(30,-10)*{}};
{\ar@/^1pc/(50,-10)*{};(45,-15)*{}};
{\ar@/^6pc/(35,-15)*{};(45,5)*{}};
\endxy
\]

\noindent The left side has $R'=4$ and $SC'=3$ and the right side has $R'=4$ and $SC' = 3$.

\noindent \textbf{Case 6:} $(54)(16)(23)$

\[
\xy
{\ar@/^-1pc/(-10,0)*{};(-5,5)*{}};
{\ar@/^-1pc/(5,5)*{};(10,0)*{}};  
{\ar@/^-1pc/(5,-15)*{};(-5,-15)*{}}; 
{\ar@/^-1pc/(-3,-4)*{};(3,-4)*{}}; 
{\ar@/^-1pc/(3,-4)*{};(-3,-4)*{}}; 
{\ar@/^-1pc/(-5,-15)*{};(5,-15)*{}};
{\ar@/^-1pc/(10,0)*{};(5,5)*{}};
{\ar@/^-1pc/(-5,5)*{};(-10,0)*{}};
(15, -4)*{\leftrightarrow};
{\ar@/^1pc/(45,5)*{};(35,5)*{}};
{\ar@/^1pc/(30,-10)*{};(35,-15)*{}};
{\ar@/^1pc/(45,-15)*{};(50,-10)*{}}; 
{\ar@/^1pc/(37,-5)*{};(43,-5)*{}}; 
{\ar@/^1pc/(43,-5)*{};(37,-5)*{}}; 
{\ar@/^-1pc/(35,5)*{};(30,-10)*{}};
{\ar@/^-1pc/(35,-15)*{};(45,-15)*{}};
{\ar@/^-1pc/(50,-10)*{};(45,5)*{}};
\endxy
\]

\noindent The left side of the diagram has $R' = 5$ and $SC' = 4$. The right side of the diagram $R'=3$ and $SC' =2$. 

\noindent \textbf{Case 7:} $(54)(16)(23)$

\[
\xy
{\ar@/^-1pc/(-10,0)*{};(-5,5)*{}};
{\ar@/^-1pc/(5,5)*{};(10,0)*{}};  
{\ar@/^-1pc/(5,-15)*{};(-5,-15)*{}}; 
{\ar@/^-1pc/(-3,-4)*{};(3,-4)*{}}; 
{\ar@/^-1pc/(3,-4)*{};(-3,-4)*{}}; 
{\ar@/^-1pc/(-5,-15)*{};(5,-15)*{}};
{\ar@/^-1pc/(-5,5)*{};(5,5)*{}};
{\ar@/^-2pc/(10,0)*{};(-10,0)*{}};
(15, -4)*{\leftrightarrow};
{\ar@/^0.5pc/(45,5)*{};(35,5)*{}};
{\ar@/^0.5pc/(35,5)*{};(45,5)*{}};
{\ar@/^1pc/(30,-10)*{};(35,-15)*{}};
{\ar@/^2.5pc/(30,-10)*{};(50,-10)*{}};
{\ar@/^1pc/(45,-15)*{};(50,-10)*{}}; 
{\ar@/^1pc/(37,-5)*{};(43,-5)*{}}; 
{\ar@/^1pc/(43,-5)*{};(37,-5)*{}}; 
{\ar@/^-1pc/(35,-15)*{};(45,-15)*{}};
\endxy
\]

\noindent The left side of the diagram has $R' = 4$ and $SC' = 3$, while the right side of the diagram has $R' = 4$ and $SC' = 3$.

In every case we have shown that the parity of $SC$ and $w(K)$ are different. Therefore, $SC - w(K) - 1$ is always even. This completes the proof of the Lemma.
\qed

%\begin{acknowledgements}
%\end{acknowledgements}


\begin{thebibliography}{}
\bibitem{amsshort}
G. Alagic, M. Jarret, and S. Jordan, Yang-Baxter operators need quantum entanglement to distinguish knots, (preprint 2015),  arXiv:1507.05979v1 21 July 2015. 

\bibitem{Alex} J. Alexander, A lemma on a system of knotted curves, Proc. Nat. Acad. Sci. USA. {\bf 9} (1923), 93-95.

\bibitem{Ara} P. K. Aravind, Borromean entanglement of the GHZ state. in
`Potentiality, Entanglement and Passion-at-a-Distance", ed. by R. S. Cohen et al,
pp. 53-59, Kluwer, 1997.

\bibitem{Baxter}  R.J. Baxter.  Exactly Solved Models in Statistical Mechanics.  Acad. Press (1982).

\bibitem{Birman} J. Birman, Braids, Links, and Mapping Class Groups, Annals of Mathematics Series Number 82,  (1974) Princeton University Press, Princeton, New Jersey.


\bibitem{Jones} 
V.F.R. Jones, A polynomial invariant for links via von Neumann algebras,
Bull. Amer. Math. Soc. {\bf 129} (1985), 103--112.

\bibitem{Jones1}
 V.F.R.Jones.  Hecke algebra representations of braid groups and link polynomials.  Ann. of Math.  126 (1987), pp. 335-338.

\bibitem{Jones2}
V.F.R.Jones.  On knot invariants related to some statistical mechanics models.  Pacific J. Math., vol. 137, no. 2 (1989), pp. 311-334.

\bibitem{Jones3} V. F. R. Jones, Braid groups, Hecke algebras and type II1 factors.
``Geometric methods in operator algebras" (Kyoto, 1983), 242Ð273, Pitman Res. Notes Math. Ser., 123, Longman Sci. Tech., Harlow, 1986.

\bibitem{Jones4}
D. Aharonov,  V. F. R. Jones, Z. Landau, A polynomial quantum algorithm for approximating the Jones polynomial. STOC'06: Proceedings of the 38th Annual ACM Symposium on Theory of Computing, 427Ð436, ACM, New York, 2006. quant-ph/0511096.

\bibitem{FKT} 
L. H. Kauffman, ``Formal Knot Theory", Princeton University Press (1983).

\bibitem{OnKnots}
L. H. Kauffman, ``On Knots", Princeton University Press (1987).

\bibitem{state}
L. H. Kauffman, State models and the Jones polynomial, \textit{Topology} \textbf{26} (1987) 395-407.

\bibitem{stat}
L. H. Kauffman, Knot theory and statistical mechanics, \textit{International Journal of Modern Physics B} \textbf{ 11}, Nos. 1, 2 (1997), 39-49.

\bibitem{virtual}
L. H. Kauffman, A self-linking invariant of virtual knots, \textit{Fundamenta Mathematicae} \textbf{184} (2004) 135-158.

\bibitem{notsoshort}
L. H. Kauffman and S. J. Lomonaco, Quantum entanglement and topological entanglement, \textit{New Journal of Physics}, 4:73.1-73.18, 2002.

\bibitem{qdeformed}
L. H. Kauffman and S. J. Lomonaco, q-Deformed Spin Networks, Knot Polynomials and Anyonic Topological Quantum Computation, \textit{Journal of Knot Theory and Its Ramifications} \textbf{16} (2007) 267-332.

\bibitem{May}
L. H. Kauffman and S. J. Lomonaco, Braiding operators are universal quantum gates,  \textit{New Journal of Physics},  6:134.1-134.40, 2004. arXiv:quant-ph/0401090.


\bibitem{tele}
L. H. Kauffman, Teleportation Topology, quant-ph/0407224, (in the Proceedings
of the  2004 Byelorus Conference on Quantum Optics),  {\it Opt. Spectrosc.} 9, 2005, 227-232.

\bibitem{quantumcompute}
L. H. Kauffman and S. J. Lomonaco, A Three-stranded quantum algorithm for the Jones polynomial, \textit{Quantum Information and Quantum Computation V}, Proceedings of Spie, April 2007, edited by E.J. Donkor, A.R. Pirich and H.E. Brandt, pp.65730T1-17, Intl Soc. Opt. Eng.

\bibitem{fibonacci}
L. H. Kauffman and S. J. Lomonaco, The Fibonacci model and the Temperley-Lieb algebra, \textit{International Journal of Modern Physics B} Vol. 22 No. 29 (2008) 5065-5080.

\bibitem{QK1} L. H. Kauffman and S. J. Lomonaco Jr., Quantum knots, in {\it Quantum Information
and Computation II, Proceedings of Spie, 12 -14 April 2004} (2004), ed. by Donkor Pirich and Brandt, pp. 268-284.

\bibitem{QK2} S. J. Lomonaco and L.  H. Kauffman, Quantum Knots and Mosaics, 
Journal of Quantum Information Processing, Vol. 7, Nos. 2-3, (2008), pp. 85 - 115. 
http://arxiv.org/abs/0805.0339

 \bibitem{QK3} S. J. Lomonaco and L. H. Kauffman, Quantum knots and lattices, or a blueprint for quantum systems that do rope tricks. Quantum information science and its contributions to mathematics, 209Ð276, Proc. Sympos. Appl. Math., 68, Amer. Math. Soc., Providence, RI, 2010.
 
 \bibitem{QK4} S. J. Lomonaco and L. H. Kauffman, Quantizing braids and other mathematical structures: the general quantization procedure. In Brandt, Donkor, Pirich, editors, {\em Quantum Information and Comnputation IX - Spie Proceedings, April 2011}, Vol. 8057, of Proceedings of Spie, pp. 805702-1 to 805702-14, SPIE 2011.

\bibitem{QK5} L. H. Kauffman and S. J. Lomonaco, Quantizing knots groups and graphs. In Brandt, Donkor, Pirich, editors, {\em Quantum Information and Comnputation IX - Spie Proceedings, April
2011}, Vol. 8057, of Proceedings of Spie, pp. 80570T-1 to 80570T-15, SPIE 2011.


\bibitem{TM} L. H. Kauffman, Knot diagrammatics, in ``Handbook of Knot Theory",
Edited by William Menasco and Morwen Thistlethwaite, pp. 233-318, (2005) Elsevier B. V. 

\bibitem{ER} Maldacena, Juan; Susskind, Leonard (2013). "Cool horizons for entangled black holes". Fortsch. Phys. 61: 781Ð811.

\bibitem{Reid} K. Reidemeister, Knotentheorie, (first published in 1932, Julius Springer, Berlin) Chelsea Publishing Company, New York (1948) 

\bibitem{RT1} 
N.Y. Reshetikhin and V. Turaev. Ribbon graphs and their invariants derived from quantum groups. Comm. Math. Phys. 127 (1990). pp. 1-26.

\bibitem{RT2}
N.Y. Reshetikhin and V. Turaev.  Invariants of Three Manifolds via link polynomials and quantum groups.  Invent. Math. 103, 547-597 (1991).

\bibitem{Yetter} D. Yetter, Markov algebras, "Braids" (Santa Cruz, CA, 1986), 705Ð730, Contemp. Math., 78, Amer. Math. Soc., Providence, RI, 1988. 

\end{thebibliography}


\end{document}